\documentclass[11pt]{amsart}

\usepackage[utf8]{inputenc}
\usepackage[T1]{fontenc}
\usepackage[english]{babel}
\usepackage[babel]{csquotes}
\usepackage{DKstyle}
\usepackage{iftex}
\ifPDFTeX
  \usepackage{bbm}
\else
  \usepackage{dsfont}
  \let\mathbbm\mathds
\fi
\usepackage[pdfstartview=FitW,colorlinks=true,linkcolor=blue,%
urlcolor=blue,citecolor=blue]{hyperref}
\usepackage{stmaryrd}
\usepackage{cleveref}
\usepackage{mathtools}
\usepackage{soul}
\usepackage{verbatim}
\usepackage{tikz-cd}
\usepackage{url}
\usepackage{xcolor}
\usepackage{enumitem}

\theoremstyle{plain}
\newtheorem{thmx}{Theorem}

\theoremstyle{definition}
\newtheorem{grassdefinition}[Thm]{Definition}

\newtheorem{fact}{Fact}[section]
\newtheorem{grassremark}[Thm]{Remark}

\newcommand{\eps}{\varepsilon}
\newcommand{\kg}{\mathfrak{g}}
\renewcommand{\dd}{\,\mathrm{d}}

\renewcommand{\be}{\mathbf{e}}

\DeclareMathOperator{\Span}{Span}
\DeclareMathOperator{\diam}{diam}

\DeclareMathOperator{\Lie}{Lie}
\DeclareMathOperator{\Var}{Var}
\DeclareMathOperator{\Cov}{Cov}

\title[A Schmidt counting theorem for Grassmannians]{A Schmidt counting theorem for Grassmannians}

\author{Ren\'e Pfitscher \and Shucheng Yu}

\address{Ren\'e Pfitscher, School of Mathematical Sciences, University of Science and Technology of China}
\email{pfitscher@edu.ustc.cn}

\address{Shucheng Yu, School of Mathematical Sciences, University of Science and Technology of China}
\email{yusc@edu.ustc.cn}

\date{}

\thanks{R.P. and S.Y. are supported by the National Key R\&D Program of China No. 2024YFA1015100.}

\begin{document}

\begin{abstract}
We prove a Schmidt counting theorem for approximation on real Grassmannians by rational subspaces of possibly different dimension. For arbitrary non-increasing approximation functions satisfying the natural divergence condition, we obtain an almost-everywhere counting asymptotic with a power-saving error term. Our approach combines homogeneous dynamics with fractional moment estimates for the Grassmann Siegel transform, overcoming the possible failure of second-moment integrability.
\end{abstract}
\maketitle

\setcounter{tocdepth}{1}
\tableofcontents

\section{Introduction}
Let $\psi : \N \rightarrow (0,\infty)$ be a non-increasing function. By Khintchine's theorem~\cite{K26}, the inequality
\begin{equation} \label{eq:Cassels}
0 \leq q x - p < \psi(q)
\end{equation}
has infinitely (resp. at most finitely) many integer solutions $(p,q) \in \Z \times \N$ for Lebesgue-almost every $x \in \R$, if the series $\sum_{q \geq 1} \psi(q)$ diverges (resp. converges). Indeed, Khintchine proved the theorem with $|q x - p| < \psi(q)$, and the one-sided form in \eqref{eq:Cassels} is due to Cassels~\cite{Cassels1950I}.

In the divergence case, Schmidt~\cite{S60} obtained a quantitative refinement. For $x\in\R$ and $T\geq1$, write
\begin{align} \label{eq:SchmidtCountingFunction}
\cN_\psi(x,T) = \# \left \{ (p,q) \in \Z \times \N : 0 \leq q x - p < \psi(q), \, 1 \leq q < T \right \}.
\end{align}
Schmidt proved that for Lebesgue-almost every $x \in \R$, as $T \rightarrow + \infty$,
$$
\cN_\psi(x,T) \, \sim \, \sum_{1\leq q < T} \psi(q).
$$
Indeed, Schmidt proved this asymptotic with an explicit error term and established analogous results in higher dimensions.

\smallskip
In \cite{S26}, de Saxc\'e studied rational approximation on real Grassmannians and established, in particular, a Khintchine-type dichotomy. In this paper, we prove a Schmidt-type counting theorem in this setting for arbitrary non-increasing approximation functions satisfying the divergence condition, with a power-saving error term. Our result allows the dimension of the subspaces being approximated to differ from that of the rational approximating subspaces, thus covering the general setting studied by Schmidt in his foundational
paper~\cite{S67}.

\subsection{Main result}\label{sec:main_result}
Let $k,\ell,d\in\N^*$ with $1\leq k,\ell<d$, and let $X_\ell=\mathrm{Gr}_{\ell,d}(\R)$ be the Grassmannian of $\ell$-dimensional linear subspaces of $\R^d$. Define $X_k$ similarly, and denote by $X_k(\Q)$ its subset of rational subspaces. In this paper, we study the approximation of $\ell$-dimensional subspaces by rational $k$-dimensional subspaces. Following Schmidt~\cite{S67}, we use the Pl\"ucker embedding $\iota_k:X_k\hookrightarrow\mathbb P(\wedge^k\R^d)$ to define a natural height function $H:X_k(\Q)\to\R_{>0}$. This height plays the role of the denominator of a rational number written in reduced form, as described in Section~\ref{sec:Height}. Let $\sigma$ and $\sigma_k$ denote the unique $\SO(d)$-invariant probability measures on $X_\ell$ and $X_k$, respectively. The measure $\sigma$ plays the role of Lebesgue measure. The quality of a rational approximation $v\in X_k(\Q)$ to $x\in X_\ell$ is measured by a natural distance $d(v,x)$ defined in terms of the principal angles between the two subspaces, see Section~\ref{sec:principal_angles}.

Since orthogonal complementation preserves both distance and height,
$$
d(v,x)=d(v^\perp,x^\perp), \qquad H(v)=H(v^\perp),
$$
we may assume throughout that $k\leq\ell$, as explained in Remark~\ref{rmk:duality}. For $x\in X_\ell$, define its \emph{Diophantine exponent} for approximation by rational $k$-dimensional linear subspaces by
$$
\beta_k(x) :=\inf\left\{ \beta>0: \exists\,c>0\ \forall\,v\in X_k(\Q),\ d(v,x)\geq cH(v)^{-\beta} \right\}.
$$
By~\cite[Theorem~1]{S26}, this exponent equals $\beta_k= \frac{d}{k(d-\ell)}$ for $\sigma$-almost every $x\in X_\ell$.

\smallskip
Analogously to~\eqref{eq:SchmidtCountingFunction}, for a function $\psi:[1,\infty)\to(0,\infty)$, $x\in X_\ell$, and $T\geq1$, define
$$
\cN_\psi(x,T) :=\#\left\{ v\in X_k(\Q): d(v,x)<\psi(H(v)),\quad H(v)<T \right\}.
$$
Let $x_0:=\Span_\R ( e_1,\ldots,e_\ell)$, where $(e_i)_{i=1}^d$ is the standard basis of $\R^d$, and set
\begin{equation}\label{eq:definition_S0}
S_0:=\{z\in X_k:z\subseteq x_0\} \simeq \mathrm{Gr}_{k,\ell}(\R).
\end{equation}
This is a Schubert variety of codimension $D:=k(d-\ell)$ in $X_k$. For non-increasing $\psi$, let
\begin{equation} \label{eq:volume_asymptotic}
M_\psi(T) := 2d\,a_{k,d}\int_1^T v(\psi(t))\,t^{d-1}\,\dd t,  
\end{equation}
where $v(r)$ is the $\sigma_k$-measure of the tube of radius $r$ around $S_0$, and $a_{k,d}>0$ is defined in~\eqref{eq:sublattice_constant_definition}. By $\SO(d)$-invariance, the tube volume is independent of the choice of $x_0\in X_\ell$. In Section~\ref{sec:Volume}, we identify $M_\psi(T)$ with the $\lambda$-volume of the counting region $\cE_\psi(T)$ defined in~\eqref{eq:Definition_E}.

\begin{thmx}\label{thm:Grassmann}
Assume $k\leq\ell$. There exists $\gamma=\gamma(k,\ell,d)\in(0,1)$ such that, for every non-increasing function $\psi:[1,+\infty)\to(0,+\infty)$ satisfying
\begin{equation}\label{eq:divergence}
\int_1^{+\infty}\psi(t)^D t^{d-1}\,\dd t=+\infty,
\end{equation}
we have, for $\sigma$-almost every $x\in X_\ell$, as $T\to\infty$,
\begin{equation}\label{eq:asymptotic}
\cN_\psi(x,T) =\tfrac12 M_\psi(T) +O_{x,\psi}\bigl(M_\psi(T)^{1-\gamma}\bigr).
\end{equation}
\end{thmx}

In the sequel, we replace $\psi$ by
$\min\{\psi,r_*\}$ where $r_* := \diam X_k + 1$ and retain the notation $\psi$. This leaves both $\cN_\psi$ and $M_\psi$ unchanged. We can make the counting formula~\eqref{eq:asymptotic} explicit by evaluating the tube volume $v(r)$. Indeed, define
$$
\Psi(T):=\int_1^T\psi(t)^D t^{d-1}\,\dd t\quad\text{and} \quad \Xi(T):=\int_1^T\psi(t)^{D+2}t^{d-1}\,\dd t,
$$
and set 
\begin{align}\label{eq:main_constant}
c_{k,\ell, d}:=\frac{\pi^{k(d-k)/2}}
     {\Gamma\!\left(1+\frac{D}2\right)}
\left(\prod_{j=1}^{k}
\frac{\Gamma(j/2)}{\Gamma((\ell-k+j)/2)}\right)
\displaystyle\prod_{j=2}^{k}\zeta(j)
     \displaystyle\prod_{j=d-k+1}^{d}\zeta(j)^{-1}.
\end{align}
By an asymptotic formula for $v(r)$ in Lemma~\ref{lem:balls} below, we will see that 
$$
M_\psi(T)\asymp \Psi(T) \quad\text{and}\quad \tfrac12M_\psi(T)=c_{k,\ell,d}\Psi(T)+O(\Xi(T)) ,
$$
where the implicit constants depend only on $k,\ell,d$. Consequently, under the divergence assumption~\eqref{eq:divergence}, Theorem~\ref{thm:Grassmann} implies $\cN_\psi(x,T)\to\infty$ for $\sigma$-almost every $x$, recovering the divergence conclusion of the Khintchine-type theorem in~\cite[Theorem~3]{S26}. Combining the second estimate with Theorem~\ref{thm:Grassmann}, we obtain
\begin{align}\label{equ:macouser}
\cN_\psi(x,T)=c_{k,\ell,d}\Psi(T) +O_{x,\psi}\bigl(\Psi(T)^{1-\gamma}\bigr) +O(\Xi(T))
\end{align}
for $\sigma$-almost every $x$. If $\psi(t)\to 0$, the term $\Xi(T)$ is negligible compared with $\Psi(T)$, so that $\cN_\psi(x,T)\sim c_{k,\ell,d}\Psi(T)$ for $\sigma$-almost every $x$.

Moreover, if $\psi(T)^2\ll_{\psi,\delta}\Psi(T)^{-\delta}$ as $T\to\infty$ for some $0<\delta<1$, then, using $\dd\Xi(t)=\psi(t)^2\,\dd\Psi(t)$, we have for any $T\geq \Psi^{-1}(1)$,
$$
\Xi(T)\ll_{\psi, \delta}1+\int_{\Psi^{-1}(1)}^T\Psi(t)^{-\delta}\,\dd \Psi(t)=1+\int_1^{\Psi(T)}u^{-\delta}\,\dd u\asymp_{\delta} \Psi(T)^{1-\delta}.
$$
Hence, for $\sigma$-almost every $x$,
$$
\cN_\psi(x,T)=c_{k,\ell,d}\Psi(T)
+O_{x,\psi,\delta}\bigl(\Psi(T)^{1-\min\{\gamma,\delta\}}\bigr).
$$
In particular, for power functions $\psi(t)=ct^{-\tau}$ with $0<c\leq r_*$ and $0<\tau\leq \beta_k$, this estimate holds with $\delta=\min\{1/2,2\tau/d\}$. In the equal-dimensional case $k=\ell$, these power-function cases were previously treated by the first-named author in~\cite[Theorem~A]{P24} for $0<\tau<\beta_k$ and~\cite[Theorem~E]{P25} for $\tau=\beta_k$, as part of more general results on rank-one flag varieties.

\subsection{Proof strategy}\label{sec:strategy}
We briefly describe the proof of Theorem~\ref{thm:Grassmann}.
Write $G=\SL_d(\R)$, $\Gamma=\SL_d(\Z)$, $\Omega=G/\Gamma$, and $K=\SO(d)$, and let $\mu$ be the $G$-invariant probability measure on $\Omega$. We first express the counting function $\cN_\psi(x,T)$ as half of a Grassmann Siegel transform on $\widetilde X_k$, evaluated on the rotation orbit $K\Gamma$, see~\eqref{equ:coufor}. Following Kelmer--Yu~\cite[\S 3,5]{KelmerYu2023a}, we then use thickening along the parabolic subgroup stabilizing $x_0$ and a Borel--Cantelli argument to deduce Theorem~\ref{thm:Grassmann} from the Haar deviation estimate in Proposition~\ref{prop:Haar}; see Section~\ref{sec:counting}.

\smallskip
The standard approach to Haar deviation estimates is to bound the variance of the Siegel transform and apply Chebyshev's inequality. Such variance bounds are typically obtained from first- and second-moment formulas; see Schmidt~\cite{Schmidt1960} and, for subsequent applications, e.g.,~\cite{AthreyaMargulis2018,AlamGhoshYu2021,H23, KelmerYu2023a,BHZ25b}. For $1<k<d-1$, however, Grassmann Siegel transforms can have infinite second moment even for bounded compactly supported functions \cite[Theorem~C]{P25}. The main new ingredient is a method for proving the required Haar deviation estimate without a second-moment formula for these transforms. We combine $L^{1+\eps}$ integrability~\cite[Theorem~A]{P25} with truncation, effective mixing, and equidistribution of certain expanding orbits.

\smallskip
More precisely, we wish to bound
$$
\mu\bigl(|\widehat{\mathbbm{1}}_{\cE} -\mu(\widehat{\mathbbm{1}}_{\cE})|\geq u \bigr),
$$
where $\cE\subseteq\widetilde X_k$ is the counting region determined by the Diophantine inequalities defining $\cN_\psi(x,T)$, and $u>0$ lies in the range specified in Proposition~\ref{prop:Haar}.
The angular condition defining $\cE$ describes a tube around $S_0$, defined in~\eqref{eq:definition_S0}, whose codimension in $X_k$ is $D=k(d-\ell)$. Following~\cite[\S 7.2]{BHZ25a} and~\cite[\S 9.1]{BHZ25b}, we approximate $\cE$ by inner and outer staircases and normalize their shells using the diagonal flow associated with $x_0$.

\smallskip
The geometric basis of this normalization is a zooming principle inspired by de Saxc\'e's work on Diophantine approximation \cite{S26,S21}, see also \cite{SH26}. In Section~\ref{sec:zooming}, we formulate this principle in coordinates transverse to $S_0$: the flow fixes $S_0$ pointwise and acts by uniform dilation in transverse unipotent coordinates. Thus $a(1/r)$ transforms a long and thin shell of height $h$ and tube radius $r$ into a region of controlled geometry and normalized height $s=hr^{D/d}$. We develop a quantitative form of this principle thus providing uniform size and boundary estimates for the normalized shells, see Lemma~\ref{lemma:shell_geometry}.

\smallskip
\begin{figure}[htbp]
\includegraphics[scale=0.5]{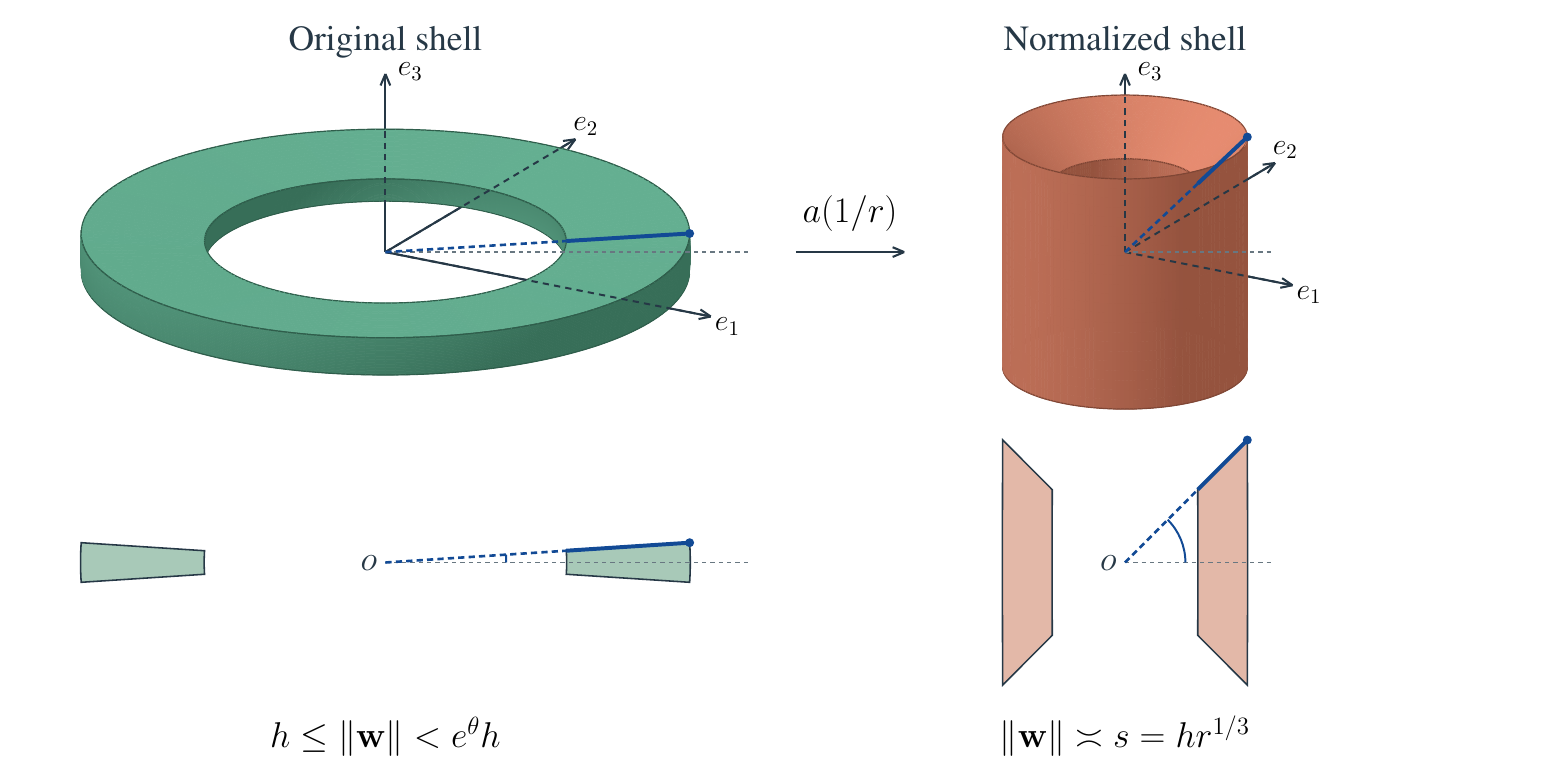}
\caption{Shell normalization for $d=3$, $\ell=2$, and $k=1$. The green shell $\mathcal D(h,r,\theta)$ on the left is one step in a staircase approximating the counting region $\cE$. Below each shell is its cross-section by the vertical plane containing the blue boundary ray. For $y>1$, the flow $a(y)=\diag(y^{-1/3},y^{-1/3},y^{2/3})$ contracts the plane $x_0=\Span(e_1,e_2)$ uniformly and expands its orthogonal complement. In projective space $X_1=\mathbb P(\R^3)$, it fixes $S_0$ pointwise and multiplies transverse unipotent coordinates by $y$ (see~\eqref{eq:zoom}), which explains the term ``zooming flow'' for $S_0$. Taking $y=1/r$ sends the green shell of small angular radius $r$ about $x_0$ to the red shell $a(1/r) \mathcal D(h,r,\theta)$ on the right, whose horizontal and vertical extents are both of order $s=hr^{1/3}$.}
\label{fig:zooming_shell}
\end{figure}

After smoothing and comparing the resulting means, it suffices to estimate
$$
\mu\bigl(|\cS-\mu(\cS)|\geq u/2\bigr),
$$
where $\cS$ is either the upper or lower smooth counting function. It has the form
$$
\cS=\sum_j\widehat f_j\circ a(e^{t_j}),
$$
where $f_j$ approximates a normalized shell of height $s_j$ and $a(e^{t_j})$ is the corresponding normalization action. We split this sum according to the normalized heights:
$$
\cS =\sum_{j\in I_{\mathrm{small}}} \widehat f_j\circ a(e^{t_j}) +\sum_{j\in I_{\mathrm{large}}} \widehat f_j\circ a(e^{t_j}).
$$
For each $j\in I_{\mathrm{small}}$, let $\Phi_j$ be a smooth truncation of $\widehat f_j$, with $0\leq\Phi_j\leq\widehat f_j$, and write
$$
\sum_{j\in I_{\mathrm{small}}} \widehat f_j\circ a(e^{t_j})= \underbrace{\sum_{j\in I_{\mathrm{small}}} \Phi_j\circ a(e^{t_j}) }_{\cS_0} + \underbrace{ \sum_{j\in I_{\mathrm{small}}} (\widehat f_j-\Phi_j)\circ a(e^{t_j}) }_{\cA}.
$$
We also define the centered large-shell contribution
$$
\cB:= \sum_{j\in I_{\mathrm{large}}} \bigl(\widehat f_j-\lambda(f_j)\bigr)\circ a(e^{t_j}).
$$
By the mean-value formula and invariance of $\mu$,
$$
\cS-\mu(\cS) =(\cS_0-\mu(\cS_0)) +(\cA-\mu(\cA))+\cB.
$$
Chebyshev's and Markov's inequalities therefore give
\begin{equation}\label{eq:intro_three}
\mu\bigl(|\cS-\mu(\cS)|\geq u/2\bigr) \ll \frac{\Var_\mu(\cS_0)}{u^2} +\frac{\|\cA\|_1+\|\cB\|_1}{u}.
\end{equation}
This reduces the problem to three estimates in Proposition~\ref{prop:three_bounds}: a variance bound for the truncated small-shell sum, an $L^1$ bound for the truncation error, and an $L^1$ bound for the large-shell discrepancy.

To bound $\Var_\mu(\cS_0)$, we apply effective mixing to the truncated functions and sum the resulting covariance bounds using the relation between normalization times and shell heights; compare \cite[proof of Proposition~7.9]{BHZ25a} and~\cite[\S 9.1]{BHZ25b}. To bound $\|\cB\|_1$, we combine equidistribution of certain expanding orbits~\cite{P24} with smoothing of the sign.

The main difficulty is to bound $\|\cA\|_1$. We dominate the Grassmann Siegel transform $\widehat f_j$ of each shell function by a smooth $K$-invariant weighted count $W_{s_j}$, itself a Grassmann Siegel transform on $\widetilde X_k$ of a smooth radial weight. The truncation at level $R\geq1$ satisfies
$$
0\leq\widehat f_j-\Phi_j \ll W_{s_j}\mathbbm{1}_{\{W_{s_j}>R\}};
$$
see~\eqref{eq:truncation_error_pointwise}. By invariance of $\mu$, this reduces the problem to bounding tail integrals of the form
$$
\int_\Omega W_s\mathbbm{1}_{\{W_s>R\}}\,\dd\mu.
$$
The desired tail estimate will follow from the following moment bound: For some fixed $\eps>0$, 
$$
\int_\Omega W_s^{1+\eps}\,\dd\mu\ll s^d\max\{1,s\}^{d\eps}\qquad(s>0),
$$
see Proposition~\ref{prop:tail}. Since a normalized shell of height $s$ and thickness $\theta$ has volume comparable to $\theta s^d$, the factor $s^d$ allows us to bound the total truncation error in terms of the sum of the shell volumes.
For the moment bound, for every $s>0$, we split $W_s$ according to the norm of the vector being counted. For vectors of norm at most one, we unfold and apply a geometric domination estimate (Lemma~\ref{lem:domination}). For the remaining vectors, we sum the annular counting bound (Lemma~\ref{lem:annular_count}) over dyadic annuli.
Combining the three estimates in~\eqref{eq:intro_three} yields Proposition~\ref{prop:Haar}.

\subsection{Notation and conventions}
We use the Landau notation $O(\cdot)$ and the Vinogradov symbol $\ll$. Given $A, B > 0$, we use the notation $A \gg B$ for $B \ll A$, and $A \asymp B$ for $A \ll B \ll A$. We use subscripts to indicate the dependence of the constant on parameters; unless stated otherwise, all implied constants are allowed to depend on $k, \ell, d$. We write $\#S$ for the cardinality of a finite set $S$.
We write $L^p(\nu)$ for the usual $L^p$-space with respect to a measure $\nu$. For $1\leq p<\infty$, we abbreviate $\|F\|_{L^p(\mu)}$ to $\|F\|_p$; norms with respect to other measures are indicated explicitly. We write $\|F\|_\infty$ for the supremum norm of $F$.

\subsection{Organization of the paper.} 
In Section~\ref{sec:preliminary}, we collect the geometric and dynamical preliminaries. In Section~\ref{sec:counting}, we deduce Theorem~\ref{thm:Grassmann} from the Haar deviation estimate in Proposition~\ref{prop:Haar}. In Section~\ref{sec:staircase}, we construct the smooth staircases and reduce the Haar estimate to the three bounds in Proposition~\ref{prop:three_bounds}. In Section~\ref{sec:tail}, we bound the tail of $W_c$ and control the truncation error. In Section~\ref{sec:shell_estimates}, we prove the small-shell variance and large-shell discrepancy bounds.

\vspace{5mm}
\textbf{Acknowledgments}.
The first-named author would like to thank Nicolas de Saxcé for proposing this problem during his PhD studies under de Saxcé’s supervision, and for his encouragement. 

\vspace{5mm}
\textbf{AI disclosure}. 
The main ideas of this paper, including the extension to arbitrary dimensions $1 \leq k, \ell < d$, were developed by the authors, and the draft was complete, except that our tail estimates for the Grassmann Siegel transform were not sufficient to be summed. A proof of the missing tail estimate in the case $\mathrm{Gr}_{2,4}(\R)$ was obtained with the help of ChatGPT (GPT-6); the domination argument underlying Proposition~\ref{prop:tail}, which extends it to all $\mathrm{Gr}_{k,d}(\R)$, was obtained with the help of Claude (Opus 5). AI was also used during the final stage of the paper for proofreading, improving the exposition, and to create Figure~\ref{fig:zooming_shell}. All statements and proofs have been checked by the authors, who take full responsibility for their correctness. 

\section{Preliminaries} \label{sec:preliminary}

\subsection{Height function} \label{sec:Height}
Let $(e_i)_{1 \leq i \leq d}$ be the standard basis of $\R^d$ and $(\mathbf e_I)_{I}$, $I\subset\llbracket 1,d\rrbracket$, $|I|=k$, the induced basis of $\wedge^k \R^d$, where $\mathbf e_I = e_{i_1}\wedge\cdots\wedge e_{i_{k}}$ for $I=\{i_1<\cdots<i_{k}\}$. The standard Euclidean structure on $\R^d$ induces a Euclidean inner product $\langle \cdot, \cdot \rangle$ on $\wedge^{k} \R^d$ for which this basis is orthonormal; denote the induced norm by $\| \cdot \|$.

The Pl\"ucker embedding $\iota_k :X_{k}\to\mathbb P(\wedge^{k}\R^d)$ sends a subspace $x$ to $[x_1\wedge\cdots\wedge x_{k}]$, where $(x_i)_{i=1}^{k}$ is any basis of $x$; we identify $X_{k}$ with its image. Recall that $x_0=\Span_\R(e_1,\ldots,e_\ell)\in X_\ell$, and let $z_0:=\Span_\R(e_1,\ldots,e_k)\in X_k$. Since $k\leq\ell$, we have $z_0\subseteq x_0$, so $z_0\in S_0$. Writing $\be_0:=e_1\wedge\cdots\wedge e_k$, we have $\iota_k(z_0)=[\be_0]$. Similarly, we identify $X_\ell$ with its image under the Pl\"ucker embedding $\iota_\ell :X_{\ell}\to\mathbb P(\wedge^{\ell}\R^d)$. 

A real linear subspace $v\subset\R^d$ is \emph{rational} if $v=\mathrm{span}_\R(v\cap\Q^d)$. For a rational ${k}$-dimensional subspace $v$, choose a $\Z$-basis $(v_i)_{i=1}^{k}$ of $v\cap\Z^d$ and define
$$
\bv:=v_1\wedge\cdots\wedge v_{k},\qquad H(v):=\|\bv\|.
$$
Note that $H(v)$ is well defined since changing the $\Z$-basis only changes $\bv$ by a sign. A nonzero integral vector $\bw\in\wedge^{k}\Z^d$ is called \emph{primitive} if $\Z\bw=\Q\bw\cap\wedge^{k}\Z^d$. The vector $\bv$ defined above is primitive and satisfies $\iota_k(v)=[\bv]$. Since $\bv$ has integer coordinates and is nonzero, we have $H(v)\geq1$. 

\subsection{The cone and the Grassmann Siegel transform} \label{sec:Siegel}
In this subsection, we define the Grassmann Siegel transform, which counts primitive pure vectors in lattices, and recall its integrability properties from \cite{P25}.

The group $G$ acts transitively on $X_{\ell}$, and the stabilizer of $x_0$ is the maximal parabolic subgroup
$$
P = \left \{ \begin{pmatrix}
A & * \\ 0 & B
\end{pmatrix} : A \in \GL_{\ell}(\R), \, B \in \GL_{d-\ell}(\R), \, \det(A) \det(B) = 1 \right \}.
$$
Thus $X_\ell\simeq P \backslash G$ via $Pg \mapsto g^{-1} x_0$. A vector of $\wedge^{k} \R^d$ is \emph{pure} if it is of the form $v_1 \wedge \dots \wedge v_{k}$ with $v_i \in \R^d$. We define the \emph{cone over $X_{k}$} by
$$
\widetilde{X}_{{k}} := \left \{ \bw \in \wedge^{k} \R^d : \bw \text{ is a nonzero pure vector} \right \} = G \cdot \be_0.
$$
We denote the set of primitive integral pure vectors by
$$
\cP_{{k}} := \left \{ \bw \in \wedge^{k} \Z^d : \bw \text{ is a primitive pure vector} \right \} = \G \cdot \be_0.
$$
The map $\cP_{k}/\{\pm 1\} \rightarrow X_{k}(\Q)$, $\bv \mapsto [\bv]$, is a bijection. The stabilizer in $G$ of the vector $\be_0$ is
$$
L = \left \{ \begin{pmatrix}
A & * \\ 0 & B
\end{pmatrix} : A \in \SL_{{k}}(\R), \, B \in \SL_{d-{k}}(\R) \right \}.
$$
The orbit map $g\mapsto g\be_0$ therefore induces a $G$-equivariant diffeomorphism
$$
G/L\longrightarrow\widetilde X_{k},\qquad gL\longmapsto g\be_0.
$$
In \cite{P25}, the first-named author extended the classical primitive Siegel transform to rational representations of semisimple $\Q$-groups; we use the following special case. Let $B_{c}^{\infty}(\widetilde{X}_{{k}})$ denote the space of bounded real-valued Borel functions with compact support on $\widetilde{X}_{{k}}$. 

\begin{definition} \label{def:Siegel}
For $f\in B_c^\infty(\widetilde X_{k})$ or a nonnegative Borel function $f$ on $\widetilde X_{k}$, define its \emph{Grassmann Siegel transform} by
$$
\widehat{f}(g\G) := \sum_{\bv \in \cP_{k}} f(g\bv) \qquad (g \in G).
$$
\end{definition}
For $f\in B_c^\infty(\widetilde X_{k})$, the sum has finitely many nonzero terms since $g\cP_{k}$ lies in the lattice $g(\wedge^{k}\Z^d)$; for $f$ a nonnegative Borel function, the sum takes values in $[0,\infty]$. In either case, it is independent of the representative $g$ since $\Gamma$ permutes $\cP_{k}$. 

We normalize the Haar measure $m_G$ on $G$ so that the induced $G$-invariant measure on $\Omega$ is a probability measure, which we denote by $\mu$.

\begin{fact} \label{fact:L1}
By \cite[Theorem~A]{P25}, using that $P_k := \Stab_G(z_0)$ is a maximal parabolic subgroup, the following hold:
\begin{enumerate}[label=(\roman*)]
\item The Grassmann Siegel transform maps $B_c^{\infty}(\widetilde{X}_{k})$ into $L^1(\mu)$.
\item There is a unique (up to scaling) $G$-invariant Radon measure $\lambda$ on $\widetilde{X}_{k}$; the Grassmann Siegel transform extends to a bounded operator $L^1(\lambda) \rightarrow L^1(\mu)$, and $\lambda$ can be normalized so that
\begin{align} \label{eq:WeilIntegrationFormula*}
\int_{\Omega} \widehat{f} \, \dd \mu = \int_{\widetilde{X}_{k}} f \, \dd \lambda =: \lambda(f) \qquad (f \in L^1(\lambda)).
\end{align}
\item The Lie group $L$ is unimodular and $\G_L = \G \cap L$ is a lattice in $L$.
\item There exists $\eps > 0$ such that the Grassmann Siegel transform maps $B_{c}^{\infty}(\widetilde{X}_{k})$ into $L^{1+\eps}(\mu)$.
\end{enumerate}
\end{fact}

By \cite[Theorem~C]{P25}, for $1 < k < d-1$ the Grassmann Siegel transform does \emph{not} map $B_c^{\infty}(\widetilde{X}_{k})$ into $L^2(\mu)$, so the parameter $\eps$ in Fact~\ref{fact:L1}~(iv) necessarily satisfies $\eps<1$ in this range. For $k \in \{1, d-1\}$ it does; this is the classical second-moment theory of the primitive Siegel transform.

\subsection{Measures, coordinates and unfolding} \label{sec:Measure}
For $y > 0$ let
$$
a(y) = \diag\left( y^{-(d-\ell)/d}I_\ell,\, y^{\ell/d}I_{d-\ell} \right), \qquad A = \{a(y) : y > 0\},
$$
and
$$
a_k(y) = \diag\left( y^{-(d-k)/d}I_k,\, y^{k/d}I_{d-k} \right), \qquad A_k := \{a_k(y) : y > 0\}.
$$
We will use $a(y)$ for the shell normalization, and $a_k(y)$ to parametrize the cone $\widetilde X_k$ in the unfolding formula below. Recall that $D = k(d-\ell)$ and set $D_k := k(d-k)$, which is equal to the dimension of $X_k$. Then $a(y) \be_0 = y^{-D/d} \be_0$ and $a_k(y) \be_0 = y^{-D_{k} /d} \be_0$. Let $m_K$ be the Haar probability measure on $K$, and $m_L$ the Haar measure on $L$ inducing a probability measure $m_{L/\G_L}$ on $L/\G_L$ where $\G_L := \Gamma \cap L$.

Under $G/L\simeq\widetilde X_k$, $gL\mapsto g\be_0$, we regard $\lambda$ as a measure on $G/L$. By the quotient integration formula~\cite[Theorem~2.51]{Folland15}, we have
\begin{equation} \label{eq:Measure_Decomposition_G}
\int_G f(g)\,\dd m_G(g)
=\int_{G/L}\int_L f(gl)\,\dd m_L(l)\,\dd\lambda(gL)
\qquad(f\in C_c(G)).
\end{equation}
By the unfolding argument below with $F\equiv1$, we will see that the quotient measure in~\eqref{eq:Measure_Decomposition_G} agrees with the normalization of $\lambda$ in~\eqref{eq:WeilIntegrationFormula*}.

Let $K_L = K \cap L = \SO(k) \times \SO(d-k)$, and let $\widetilde{\sigma}_k$ be the push-forward of $m_K$ to $K/K_L$, the Grassmannian of oriented $k$-dimensional linear subspaces; the projection $K/K_L \rightarrow X_{k}$ pushes $\widetilde{\sigma}_k$ forward to the unique $K$-invariant probability measure $\sigma_k$ on $X_k$. By \cite[\S 2.7, Equation~(2.13)]{P24}, the map $(\mathrm{k} K_L, a_k(y)) \mapsto \mathrm{k} a_k(y) \be_0$ is a homeomorphism $(K/K_L) \times A_k \rightarrow \widetilde{X}_{k}$, and for some $\omega_{k,d}>0$ the cone measure is given by
\begin{equation} \label{eq:Measure-On_Tilde-X}
\dd \lambda(\mathrm{k} a_k(y) \be_0) = \omega_{k,d} \, y^{-(D_{k}  +1)} \, \dd \widetilde{\sigma}_k(\mathrm{k}) \, \dd y .
\end{equation}
For $s > 0$ let $B_s := \{\bw \in \widetilde{X}_k : \|\bw\| \leq s\}$. 
By the mean value formula~\eqref{eq:WeilIntegrationFormula*}, the measure description~\eqref{eq:Measure-On_Tilde-X} and the relation $\|\mathrm{k} a_k(y)\be_0\|=y^{-D_{k} /d}$, we have
$$
\int_\Omega\widehat{\mathbbm{1}}_{B_s}\,\dd\mu=\lambda(B_s)=\omega_{k,d}\int_{s^{-d/D_{k} }}^\infty y^{-D_{k} -1}\,\dd y
=\frac{\omega_{k,d}}{D_{k} } s^d.
$$
On the other hand, the mean value formula for primitive sublattices in \cite[Theorem~3]{K22} implies
$\int_{\Omega} \widehat{\mathbbm{1}}_{B_s}\, \dd \mu = 2 \, a_{k,d} \, s^d$,
with
\begin{equation}\label{eq:sublattice_constant_definition}
a_{k,d}:= \frac{1}{d} \binom{d}{k} \left ( \prod_{i = 1}^k \frac{V_{d-i+1}}{V_i} \right ) \prod_{j = 2}^{k} \zeta(j) \prod_{j = d - k +1}^{d} \zeta(j)^{-1}.
\end{equation}
Here $V_j:=\pi^{j/2}/\Gamma(1+j/2)$ denotes the volume of the Euclidean unit ball in $\R^j$. Note that the factor $2$ here accounts for the two primitive pure vectors $\pm \bv$ representing a primitive sublattice. Hence
\begin{equation} \label{eq:sublattice_constant}
\lambda(B_s) = 2 a_{k,d} \, s^d \quad \text{and}\quad \omega_{k,d} = 2 D_{k}  \, a_{k,d}.
\end{equation}

By~\eqref{eq:Measure-On_Tilde-X} and~\eqref{eq:sublattice_constant}, making the change of variables $r=y^{-D_{k} /d}$ and projecting $\widetilde{\sigma}_k$ to $\sigma_k$, we obtain, for every nonnegative Borel function $F$ on $(0,\infty)\times X_k$,
\begin{equation}\label{eq:polar_measure}
\int_{\widetilde X_k}F(\|\bw\|,[\bw])\,\dd\lambda(\bw)
=2da_{k,d}\int_0^\infty\int_{X_k}F(r,z)\,\dd\sigma_k(z)\,r^{d-1}\dd r.
\end{equation}

We now record an unfolding identity for the Grassmann Siegel transform, which will be used to prove the tail and large-shell estimates in Sections~\ref{sec:tail} and~\ref{sec:large}.

\begin{lemma} \label{lem:unfolding}
For all non-negative Borel functions $\varphi$ on $\widetilde{X}_{k}$ and $F$ on $\Omega$,
{
\begin{multline*}
\int_{\Omega} \sum_{\bv \in \cP_{k}} \varphi(g\bv) \, F(g\G) \, \dd \mu(g\G) \\
= \omega_{k,d} \int_{K/K_L} \int_0^{\infty} \varphi \big ( \mathrm{k} a_k(y) \be_0 \big ) \left ( \int_{L/\G_L} F \big ( \mathrm{k} a_k(y) l \G \big ) \, \dd m_{L/\G_L}(l) \right ) \frac{\dd y}{y^{D_{k} +1}} \, \dd\widetilde{\sigma}_k(\mathrm{k}K_L).
\end{multline*}
}
In particular, taking $F\equiv1$ recovers the mean value formula~\eqref{eq:WeilIntegrationFormula*}.
\end{lemma}
\begin{proof}
Recall that $m_G$ locally agrees with $\mu$. We also write $m_G$ for its induced measure on $G/\Gamma_L$, and regard $\lambda$ as the quotient measure in~\eqref{eq:Measure_Decomposition_G}. Since $\gamma\Gamma_L\mapsto\gamma\be_0$ identifies $\Gamma/\Gamma_L$ with $\cP_k$, unfolding gives
\begin{align*}
\int_\Omega\widehat\varphi(g\Gamma)F(g\Gamma)\,\dd\mu
&=\int_{G/\Gamma_L}\varphi(g\be_0)F(g\Gamma)\,\dd m_G\\
&=\int_{G/L}\varphi(g\be_0)
\left(\int_{L/\Gamma_L}F(gl\Gamma)\,\dd m_{L/\Gamma_L}(l)\right)\dd\lambda(gL).
\end{align*}
Here the second equality follows from~\eqref{eq:Measure_Decomposition_G} after passing to the quotient by $\Gamma_L$, using $l\be_0=\be_0$.
Applying~\eqref{eq:Measure-On_Tilde-X} to the quotient integral gives the stated formula.

Taking $F\equiv1$ and using $m_{L/\Gamma_L}(L/\Gamma_L)=1$, we obtain the mean value formula for bounded compactly supported functions on $\widetilde X_k$. We can then extend it to all $L^1$-functions by a standard approximation argument.
\end{proof}

\subsection{Principal angles and distance} 
\label{sec:principal_angles}
In this section, we specify the distance $d(\cdot, \cdot)$ that we use to measure the quality of a rational approximation $v \in X_k(\Q)$ to a point $x \in X_\ell$. We shall always identify the points of the Grassmannian with the corresponding subspaces of $\R^d$.

\smallskip
Let $x,z \subseteq \R^d$ be subspaces of dimensions $1 \leq k, \ell < d$ and let $m := \min \{k,\ell\}$. The relative position of $x$ and $z$ is described by their \emph{principal angles} (see \cite{J75}), defined geometrically as follows. The first principal angle $\theta_1 \in [0, \frac{\pi}{2}]$ is the smallest angle between a unit vector of $x$ and a unit vector of $z$, and it is attained by a pair $(a_1, b_1)$. Recursively, $\theta_i$ is the smallest such angle among unit vectors $a \in x$ and $b \in z$ with $a \perp a_1, \dots, a_{i-1}$ and $b \perp b_1, \dots, b_{i-1}$, attained by a pair $(a_i, b_i)$. This yields
$$
0 \leq \theta_1 \leq \dots \leq \theta_m \leq \tfrac{\pi}{2},
$$
together with orthonormal \emph{principal vectors} $a_1, \dots, a_m$ of $x$ and $b_1, \dots, b_m$ of $z$. Thus $\theta_1 = 0$ if and only if $x \cap z \neq \{0\}$, and $\theta_m = 0$ if and only if $x \subset z$ or $z \subset x$. We denote by $\theta(x,z) = (\theta_1, \dots, \theta_m)$ the associated vector of principal angles.

\smallskip
We then define the distance between $x$ and $z$ to be the Euclidean norm of $\theta(x,z)$:
\begin{equation} \label{eq:distance_definition}
d(x,z) := \| \theta(x,z) \|_2 .
\end{equation}

\begin{remark} \label{rmk:duality}
The map $v \mapsto v^\perp$ defines a bijection $X_k(\Q)\rightarrow X_{d-k}(\Q)$ preserving height,
$$
H(v^\perp)=H(v),
$$
see \cite[Equation (4), p. 433]{S67}. Moreover, by \cite[Theorem~2.7, Equation~(2.3)]{KA06}, we have
$$
d(x,z) = d(x^\perp, z^\perp), \qquad \text{for all subspaces $x,z \subseteq \R^d$.}
$$
Writing $\cN_\psi^{k,\ell}$ for the counting function for approximation of $\ell$-subspaces by rational $k$-subspaces, we then have
$$
\cN_\psi^{k,\ell}(x,T)
=\cN_\psi^{d-k,d-\ell}(x^\perp,T)
\qquad(x\in X_\ell,\ T\geq1).
$$
If $k>\ell$, then $d-k<d-\ell$. Since orthogonal complementation carries the invariant probability measure on $X_\ell$ to that on $X_{d-\ell}$, almost-everywhere statements are preserved under orthogonal complementation. The counting formula in Theorem~\ref{thm:Grassmann} therefore extends to the case $k>\ell$ via this identity, by applying the theorem to the complementary dimensions. 
\end{remark}

From now on, we assume $k\leq\ell$, as justified in Remark~\ref{rmk:duality}. We now collect some facts about this distance, along with precise references.

\begin{fact} \label{fact:principal_angles}
The following hold:
\begin{enumerate}[label=(\roman*)]
\item When restricted to $X_k\times X_k$, the distance $d$ is induced by the $K$-invariant Riemannian metric obtained by pulling back, under the Pl\"ucker embedding $\iota_k$, the standard projective Riemannian metric induced by the unit sphere equipped with the metric inherited from the Euclidean inner product; see \cite[\S 2, pp.~1179--1180]{YL16}. We denote the associated Riemannian volume measure by
$\mathrm{vol}$. By uniqueness of the $K$-invariant probability
measure on $X_k$, we have $\sigma_k = \vol(X_k)^{-1} \vol$.
\item The total volume of $X_k$ is given by
\begin{equation}\label{eq:Volume_Grassmannian}
\vol(X_k) =\pi^{D_{k} /2}\prod_{j=1}^k\frac{\Gamma(j/2)}{\Gamma((d-k+j)/2)} =\binom dk\prod_{i=1}^k\frac{V_{d-i+1}}{V_i} ,
\end{equation}
where $V_j = \pi^{j/2}/\Gamma(1+j/2)$ denotes again the volume of the Euclidean unit ball in $\R^j$, see~\cite[\S 9, Equation~(23)]{YL16}. 
\item For every $x \in X_\ell$, the set $S(x) := \{z \in X_k : z \subseteq x\}$ is a Schubert variety of codimension $k(d-\ell)$, see~\cite[\S 8]{YL16}. Moreover, for every $z \in X_k$ and $x \in X_\ell$,
\begin{equation} \label{eq:distance_Schubert}
    d(z,x) = \min_{z' \in S(x)} d(z,z') =: d(z,S(x)),
\end{equation}
see~\cite[\S 4, Equations~(11), (14)]{YL16}. In the current notation, $S_0=S(x_0)$, where $S_0$ is defined in~\eqref{eq:definition_S0}.
\item We equip $S_0$ with the restriction of the Riemannian metric on $X_k$, and denote its Riemannian volume measure by $\vol_{S_0}$. By \cite[Propositions~21, 22]{YL16}, $S_0$ is isometrically identified with $\mathrm{Gr}_{k,\ell}(\R)$. The total volume is thus
\begin{equation}\label{eq:Volume_Schubert}
\vol_{S_0}(S_0) =\pi^{k(\ell-k) /2}\prod_{j=1}^k\frac{\Gamma(j/2)}{\Gamma((\ell-k+j)/2)} =\binom \ell k \prod_{i=1}^k\frac{V_{\ell-i+1}}{V_i}.
\end{equation}
\end{enumerate}
\end{fact}

\smallskip
For every $r > 0$ we denote
$$
v(r):=\sigma_k \left (\{z\in X_k:d(z,S_0)<r\} \right) 
$$
the normalized volume of the tube of radius $r$ around $S_0$. The following basic lemma gives the small-radius volume asymptotic needed to evaluate the mean $M_\psi(T)$.
\begin{lemma}\label{lem:balls}

For every $r>0$, we have
\begin{align}\label{equ:vrcomp}
v(r)=\frac{\vol_{S_0}(S_0)V_D}{\vol(X_k)}\,r^D
+O_{d}(r^{D+2})
\quad \text{and}\quad  v(r)\ll_{d}r^D.
\end{align}
Moreover, for every $r>0$ and $0<\delta<1/2$, we have
\begin{equation}\label{eq:RiemannianVolumeDelta}
v((1\pm\delta)r)=v(r)(1+O_{d}(\delta)) .
\end{equation}
\end{lemma}

{
\begin{remark}
Combining \eqref{equ:vrcomp} with \eqref{eq:volume_asymptotic}, the leading constant $c_{k,\ell,d} > 0$ in \eqref{equ:macouser} is given by
\begin{equation} \label{eq:value_explicit}
c_{k,\ell,d} = d \, a_{k,d} \frac{\vol_{S_0}(S_0) V_D}{\vol(X_k)}.
\end{equation}
Inserting the explicit values for $a_{k,d}$ in \eqref{eq:sublattice_constant_definition} and for the volumes $\vol_{S_0}(S_0)$ and $\vol(X_k)$ in Fact~\ref{fact:principal_angles} we get the explicit formula for $c_{k,\ell,d}$ in~\eqref{eq:main_constant}.
\end{remark}
}

\begin{proof}[Proof of Lemma~\ref{lem:balls}]
Consider the rank-$D$ normal bundle
$$
NS_0 = \{(y,u):y\in S_0,\ u\in (T_yS_0)^\perp \subset T_y X_k\}.
$$
Its normal exponential map is $E : NS_0 \longrightarrow X_k$, $E(y,u)=\exp_y(u)$. Since $X_k$ is compact, it is geodesically complete by Hopf--Rinow, so $E$ is defined on all of $NS_0$. By the tubular neighborhood theorem and compactness of $S_0$, there exists $r_0>0$ such that $E$ restricts to a diffeomorphism
$$
\{(y,u)\in NS_0 :\|u\|<r_0\}
\longrightarrow
\{z\in X_k : d(z,S_0) < r_0\},
$$
such that $d(E(y,u),S_0)=\|u\|$ for all $u$ with $\|u\| < r_0$, see \cite[Theorems~5.25 and~6.40]{L18}. Here the norm $\|\cdot\|$ corresponds to the inner product on each normal fiber induced by the Riemannian metric on $X_k$.
Equip $N S_0$ with the smooth density obtained by choosing the induced Riemannian volume on $S_0$ and the Euclidean volume on each normal fiber:
$$
\int_{N S_0}f\,\mathrm{d}\nu_N
:=
\int_{S_0} \int_{N_y S_0} f(y,u) \,\mathrm{d}u\,\mathrm{d}\vol_{S_0}(y).
$$
Let $J(y,u)=|\det DE_{(y,u)}|$ be the absolute volume Jacobian of $E$ relative to $\dd \nu_N$ and the Riemannian volume measure on $X_k$. 
Thus, for every Borel set $A\subset\{(y,u)\in NS_0:\|u\|<r_0\}$, by the change-of-variables formula we have
{
$$
\vol(E(A))=\int_A J(y,u)\,\dd\nu_N(y,u).
$$
}
{
Writing $s(y)=(y,0)$ for the zero section, there is a canonical identification $T_{(y,0)}NS_0 \simeq T_yS_0\oplus N_yS_0$ via the following map: 
$$
T_yS_0\oplus N_yS_0\longrightarrow T_{(y,0)}NS_0,
\qquad (\xi,\eta)\longmapsto Ds_y(\xi)+\left.\frac{\mathrm{d}}{\mathrm{d}t}(y,t\eta)\right|_{t=0}.
$$
}
Under this identification, using also the identities $E(y,0)=y$ and $\left.\frac{\mathrm{d}}{\mathrm{d}t}E(y,t\eta)\right|_{t=0}=\eta$ we get
$$
DE_{(y,0)}(\xi,\eta)=\xi+\eta.
$$
Since this decomposition is orthogonal, this gives $J(y,0)=1$. Smoothness near the zero section and compactness of $S_0$ therefore give the uniform fiberwise Taylor expansion
$$
J(y,u)=1+L_y(u)+O(\|u\|^2) ,
$$
where $L_y$ is linear on $N_y S_0$. For $0<r<r_0$, integrating over centered normal balls cancels the linear term, yielding
$$
\begin{aligned}
v(r) =\frac1{\vol(X_k)} \int_{S_0}\int_{\|u\|<r}J(y,u)\,\mathrm{d}u\,\mathrm{d}\vol_{S_0}(y)
=\frac{\vol_{S_0}(S_0)V_D}{\vol(X_k)}r^D+O(r^{D+2}).
\end{aligned}
$$
This proves both estimates in \eqref{equ:vrcomp} for all $0<r<r_0$; we can then extend them to all $r>0$ using $v(r)\leq 1$. 

Let us now show~\eqref{eq:RiemannianVolumeDelta}. Set $R := \operatorname{diam}(X_k)$. For $0< a <b \le 3 R$, every $z$ satisfying $a\le d(z,S_0)<b$ has a nearest point $y\in S_0$, because $S_0$ is compact. A minimizing geodesic from $y$ to $z$ meets $S_0$ orthogonally at $y$, by the first variation formula. Consequently, there exists some $u\in (T_yS_0)^{\perp}$ with $z=\exp_y(u)$ and $\|u\|=d(z,y)=d(z,S_0)$. Hence
$$
\{z:a\le d(z,S_0)<b\} \subseteq E\bigl(\{(y,u)\in N S_0:a\le\|u\|<b\}\bigr).
$$
The absolute Jacobian $J$ is continuous and bounded on the compact disk bundle
$$
\{(y,u)\in N S_0:\|u\|\le 3 R \}.
$$
By the area formula with multiplicity~\cite[Theorem~C.2.1]{Howard94}, we have
$$
\begin{aligned}
v(b)-v(a) \le \frac1{\vol(X_k)} \int_{S_0}\int_{a\le\|u\|<b}J(y,u)\,\mathrm{d}u\,\mathrm{d}\vol_{S_0}(y) \ll b^D-a^D.
\end{aligned}
$$
Moreover, the small-radius asymptotic and monotonicity imply $v(r)\gg r^D$ for all $0<r\le 2 R$. It follows that for $0<r\le 2 R$ and $0<\delta<1/2$,
$$
\begin{aligned}
|v((1\pm\delta)r)-v(r)| \ll \bigl|(1\pm\delta)^D-1\bigr|r^D \ll \delta r^D \ll \delta v(r).
\end{aligned}
$$
If $r>2R$, then $(1-\delta)r> R$, so all three volumes are $1$, and the estimate holds trivially. 
\end{proof}

For $g\in G$, write $\|g\|_{\mathrm{op}}$ for its operator norm with respect to the standard Euclidean norm on $\R^d$, and $\|\wedge^k g\|_{\mathrm{op}}$ for the operator norm of the $g$-action on $\wedge^{k}\R^d$ with respect to the induced Euclidean norm on $\wedge^{k}\R^d$. Note that if $s_1\geq\cdots\geq s_d>0$ are the singular values of $g$, then
\begin{equation}\label{eq:exterior_operator_norm}
\|\wedge^k g\|_{\mathrm{op}}=s_1\cdots s_k \leq s_1^k=\|g\|_{\mathrm{op}}^k.
\end{equation}

We will use the following estimate in Lemma~\ref{lem:slice_comparison} below to compare counting regions under small perturbations in the thickening argument.

\begin{lemma}\label{lem:distance_distortion}
For all $g\in G$ and all subspaces $x,z \subseteq \R^d$ with $k = \dim z \leq \ell = \dim x$, we have
\begin{equation}\label{eq:grassmann_lipschitz_bound}
d(gz,gx)\leq\|\wedge^k g\|_{\mathrm{op}}\|\wedge^k g^{-1}\|_{\mathrm{op}}\,d(z,x).
\end{equation}
\end{lemma}

\begin{proof}
Set $A = \wedge^k g$. The projective action induced by $A$ lifts to the map on the unit sphere $F_A(v)= Av / \|Av\|$, where $v$ is a unit vector in $\wedge^k\R^d$. For a tangent vector $w$ at $v$, differentiation gives
$$
D_vF_A(w) = \frac{(I-F_A(v)F_A(v)^T)Aw}{\|Av\|}.
$$
Here $I-F_A(v)F_A(v)^T$ is the orthogonal projection onto the orthogonal complement of $F_A(v)$. Since the orthogonal projection has operator norm at most $1$, and $\|Av\|\ge\|A^{-1}\|_{\mathrm{op}}^{-1}$, we get
$$
\|D_vF_A(w)\|\leq\|A\|_{\mathrm{op}}\|A^{-1}\|_{\mathrm{op}}\|w\|.
$$
The covering map from the unit sphere in $\wedge^k\R^d$ to $\mathbb{P}(\wedge^k\R^d)$ is a local isometry, and $A$ preserves the Plücker image of $X_k$. Thus, by Fact~\ref{fact:principal_angles}~(i), the action of $g$ increases lengths of curves in $X_k$ by at most $\|A\|_{\mathrm{op}}\|A^{-1}\|_{\mathrm{op}}$. Taking the infimum over curves gives $d(gz,gz') \leq \|A\|_{\mathrm{op}}\|A^{-1}\|_{\mathrm{op}}\,d(z,z')$ for all $z,z' \in X_k$. For $x \in X_\ell$, we have $gS(x)=S(gx)$ and using the identity $d(z,x)=d(z, S(x))$ in \eqref{eq:distance_Schubert} yields
\begin{displaymath}
d(gz,gx) = \inf_{z'\in S(x)} d(gz,gz') \leq \|A\|_{\mathrm{op}}\|A^{-1}\|_{\mathrm{op}} \inf_{z'\in S(x)} d(z,z') =\|A\|_{\mathrm{op}}\|A^{-1}\|_{\mathrm{op}}\,d(z,x). \qedhere
\end{displaymath}
\end{proof}

{
For later use, we also fix a right-invariant Riemannian metric $d_G$ on $G$. For $r>0$, let
$$
B_G(r):=\{g\in G:d_G(I_d,g)<r\}
$$
be the open ball of radius $r$ centered at the identity $I_d$. In what follows, constants and auxiliary parameters may also depend on this fixed metric.
}

\subsection{Zooming} \label{sec:zooming}
The shell normalization in Section~\ref{sec:staircase} is based on a simple geometric observation: the flow $a(y)$ fixes $S_0$ pointwise and dilates the transverse directions to $S_0$ by the same factor $y$. Thus $a(1/r)$ expands a thin tube of radius $r$ around $S_0$ to a fixed transverse scale. In this subsection, we introduce suitable coordinates and express this dilation property through the zooming identity~\eqref{eq:zoom} below.

\smallskip
 Recall that $x_0=\Span_\R(e_1,\ldots,e_\ell)$ and $z_0=\Span_\R(e_1,\ldots,e_k)\subseteq x_0$. The Lie algebra $\mathfrak{u}^-$ of the unipotent subgroup $U^-$ opposite to $P=\Stab_G(x_0)$ is
$$
\mathfrak{u}^- = \left \{ u_B := \begin{pmatrix}
0 & 0 \\ B & 0
\end{pmatrix} : B \in \mathrm{M}_{(d-\ell) \times \ell}(\R) \right \} .
$$
We equip $\mathfrak{u}^-$ with the Frobenius norm $\|\cdot\|_{\mathrm F}$. The flow $a(y)$ acts on $\mathfrak{u}^-$ via the adjoint representation simply by multiplication by $y$, that is,
\begin{equation} \label{eq:adjoint_simple}
\Ad(a(y)) u = y  u , \qquad ( y > 0, \, u \in \mathfrak{u}^- ) .
\end{equation}
For coordinates transverse to $S_0$ at $z_0$, we use
$$
\mathfrak{n} := \left\{ u_{(B\;0)}:B\in\mathrm{M}_{(d-\ell)\times k}(\R)
\right\}
\simeq\mathrm{Hom}(z_0,x_0^\perp).
$$
Since $u^2=0$, we have $\exp(u)=I_d+u$. For $u=u_{(B\;0)}$, let $\mathcal L_B:z_0\to x_0^\perp$ be the linear map represented by $B$ in the standard bases. Then $\exp(u)z_0$ is the graph $\{x+\mathcal L_B(x):x\in z_0\}$. Write $\pi_{x_0} : \R^d \rightarrow \R^d$ for orthogonal projection onto $x_0$ and let
$$
\mathcal O:=\{z\in X_k:z\cap x_0^\perp=\{0\}\}.
$$
This is an open dense subset of $X_k$ containing $S_0$. For any $k$-subspace $z$, write $\pi_{x_0}(z):=\{\pi_{x_0}(v):v\in z\}\subset x_0$ for its image. This image is a $k$-subspace exactly when $\pi_{x_0}|_z$ is injective, or equivalently, when $z\in\mathcal O$. For each $z\in\mathcal O$, set $W:=\pi_{x_0}(z)$. Then $\pi_{x_0}|_z:z\to W$ is an isomorphism. Define
$$
\mathcal L_z(x):=(\pi_{x_0}|_z)^{-1}(x)-x\in x_0^\perp,\qquad x\in W.
$$
This is a linear map and $z=\{x+\mathcal L_z(x):x\in W\}$. It is the unique map with this graph, so every $z\in\mathcal O$ is uniquely the graph of a linear map $\pi_{x_0}(z)\to x_0^\perp$.
In particular, if $\pi_{x_0}(z)=z_0$, then $W=z_0$ and the graph map $\mathcal L_z$ is represented by a unique matrix $B\in\mathrm M_{(d-\ell)\times k}(\R)$. Its graph is $\exp(u_{(B\;0)})z_0$, and conversely every such graph projects onto $z_0$. Hence
\begin{equation}\label{eq:graph_projection}
\pi_{x_0}(z)=z_0 \quad\Longleftrightarrow\quad z=\exp(u)z_0 \text{ for a unique }u\in \mathfrak{n} .
\end{equation}

Set $H=\Stab_K(x_0)$, so that $S_0=Hz_0$. Since $H$ acts transitively on the $k$-planes contained in $x_0$, for each $z\in\mathcal O$ choose $h\in H$ such that $hz_0=\pi_{x_0}(z)$. Since $h$ preserves the decomposition $x_0\oplus x_0^{\perp}$, it commutes with $\pi_{x_0}$. Hence $z_0=\pi_{x_0}(h^{-1}z)$, and applying \eqref{eq:graph_projection} to $h^{-1}z$ gives some $u\in \mathfrak{n}$ with
\begin{equation}\label{eq:zooming_coordinates}
z=h\exp(u)z_0,\qquad h\in H,\quad u\in \mathfrak{n} .
\end{equation} 
Every $h\in H$ commutes with $a(y)$, and $a(y)$ fixes $z_0$. Therefore, for $z=h\exp(u)z_0$ as in \eqref{eq:zooming_coordinates}, \eqref{eq:adjoint_simple} gives the following zooming identity
\begin{equation}\label{eq:zoom}
a(y)h\exp(u)z_0=h\exp(yu)z_0 .
\end{equation}
Thus, in the representation $z=h\exp(u)z_0$, the point $hz_0\in S_0$ remains fixed under $a(y)$, while the transverse parameter $u$ is multiplied by $y$. To relate this coordinate to distance, for $u=u_{(B\;0)}$ with $B\in \mathrm{M}_{(d-\ell)\times k}(\R)$, let $s_1(B)\geq\cdots\geq s_k(B)\geq0$ be the singular values of $B$, namely the nonnegative square roots of the eigenvalues of $B^{\mathsf T}B$, counted with multiplicity. By the singular-value characterization of principal angles \cite[\S 2, equations~(4)-(5)]{YL16}, the principal angles between $\exp(u)z_0$ and $x_0$ are $\arctan s_k(B),\ldots,\arctan s_1(B)$ in non-decreasing order. Using also the facts that $h\in H\subset K$ acts isometrically on $X_k$ and preserves $S_0$ we have
\begin{align}\label{eq:paidu}
d(h\exp(u)z_0,S_0)^2 =\sum_{i=1}^k\arctan^2 s_i(B)\quad \text{and} \quad \|u\|_{\mathrm F}^2=\sum_{i=1}^k s_i(B)^2.
\end{align}
In particular, uniformly in $h\in H$,
\begin{equation} \label{eq:zooming_projective}
d(h\exp(u)z_0,S_0)\asymp\|u\|_{\mathrm F} \qquad(\|u\|_{\mathrm F} \leq 1).
\end{equation}
For sufficiently small $r>0$, the tube of radius $r$ around $S_0$ is contained in $\mathcal O$. Hence points in this tube have transverse coordinates of size $O(r)$, which $a(1/r)$ sends to a fixed bounded set. Finally, since $a(y)\be_0=y^{-D/d}\be_0$, the corresponding identity on the cone is
\begin{align}\label{eq:zooming_cone}
a(y)h\exp(u)\be_0 =y^{-D/d}h\exp(yu)\be_0.
\end{align}

\subsection{Dynamical inputs} \label{sec:Equidistribution}
Let $\kg=\mathfrak{sl}_d(\R)$ be the Lie algebra of $G$. For $Z\in\kg$ and $\phi\in C^\infty(\Omega)$, define the Lie derivative $Z\phi$ by $(Z\phi)(x):=\left.\frac{\dd}{\dd t}\phi(\exp(tZ)x)\right|_{t=0}$. Recall that  $K = \SO(d)$. Fix an orthonormal basis of $\kg$ containing an orthonormal basis of the Lie algebra $\mathfrak{k} = \Lie(K)$, and let $\cM_r$, resp.\ $\cM_r^K$, be the set of ordered monomials of degree at most $r$ in this basis, resp.\ in its elements from $\mathfrak{k}$.

For a full-rank lattice $\Lambda\subset\R^d$, define its first minimum by
\begin{align}\label{eq:first_minimum}
\lambda_1(\Lambda):=\min_{\bv\in\Lambda\setminus\{0\}}\|\bv\|,
\end{align}
where $\|\cdot\|$ is the Euclidean norm.
For $x = g\G \in \Omega$ let $\mathrm{ht}(x) := \lambda_1(g\Z^d)^{-1}$.
We use the following Sobolev norms
$$
\cS_{2,r}(\phi) := \Big ( \sum_{\mathscr{D} \in \cM_r} \big \| \mathrm{ht}^{\, r} \, \mathscr{D} \phi \big \|_{2}^2 \Big )^{1/2} 
\quad\text{and}\quad
\cS^K_{r}(\phi) := \Big ( \sum_{\mathscr{D} \in \cM^K_r} \| \mathscr{D} \phi \|_{2}^2 \Big )^{1/2}.
$$
The first is weighted in the cusp and involves all directions; the second involves only the compact directions and is unweighted. The proof uses two dynamical inputs, both consequences of the spectral gap of $G$ acting on $L^2_0(\mu)$.

\begin{fact}[Equidistribution of expanding $L$-orbits] \label{fact:equidistruibution_L}
By \cite[Lemma~5.1]{P24}, there are $c_L > 0$ and an integer $r_L \geq 1$ such that for all $\phi \in C_c^\infty(\Omega)$ and $y \in (0,1]$,
$$
\int_{L / \G_L} \phi(a_k(y) l \G) \, \dd m_{L / \G_L} = \int_{\Omega} \phi \, \dd \mu + O \big ( y^{c_L} \, \cS_{2,r_L}(\phi) \big ).
$$
\end{fact}

\begin{remark}
In \cite[\S2.8]{P24}, the Sobolev weight is defined using $\mathrm{ht}_{\mathrm{Ad}}(g\Gamma):=\lambda_1(\Ad(g)\kg_\Gamma)^{-1}$, where $\kg_\Gamma$ is a fixed rational $\Ad(\Gamma)$-invariant lattice in $\mathfrak{sl}_d(\R)$. Since $\kg_\Gamma$ is commensurable with $\mathfrak{sl}_d(\Z)$ and the fixed Euclidean norm on $\kg$ is equivalent to the operator norm, we may compute this height using $\mathfrak{sl}_d(\Z)$ and the operator norm up to a factor independent of $g$.
Let $\lambda_1\leq \cdots\leq \lambda_d$ be the successive minima of $\Lambda=g\Z^d$. Every nonzero $A\in g\mathfrak{sl}_d(\Z)g^{-1}$ maps $\Lambda$ into itself and is nonzero on at least one of $d$ independent lattice vectors of norm at most $\lambda_d$. Write $\bv$ for such a lattice vector. Then $\lambda_1\leq \|A\bv\|\leq \|A\|_{\mathrm{op}}\|\bv\|\leq  \|A\|_{\mathrm{op}}\lambda_d$. Thus $\|A\|_{\mathrm{op}}\geq\lambda_1/\lambda_d$. By Minkowski's second theorem, $\lambda_1^{d-1}\lambda_d\leq\lambda_1\cdots\lambda_d\ll_d1$, so $\mathrm{ht}_{\mathrm{Ad}}(g\Gamma)\ll\lambda_d/\lambda_1\ll \lambda_1^{-d}=\mathrm{ht}(g\Gamma)^d$. Thus, after increasing $r_L$ if necessary, $\cS_{2,r_L}$ dominates the Sobolev norm used in \cite[Lemma~5.1]{P24}. 
\end{remark}

\begin{fact}[Decay of matrix coefficients] \label{fact:mixing}
By \cite[Corollary~3.5]{KleinbockMargulis1999} (see also \cite[Corollary~2.4.4]{KleinbockMargulis1996}), there are $\kappa_{\mathrm{mix}}>0$ and an integer $m\geq1$ such that for all $\phi_1,\phi_2\in C^\infty(\Omega)$ with $\cS_m^K(\phi_i)<\infty$ and all $t\in\R$,
$$
\left | \int_{\Omega} \phi_1 \big ( a(e^t) x \big ) \, \phi_2(x) \, \dd \mu(x) - \int_{\Omega} \phi_1 \, \dd \mu \int_{\Omega} \phi_2 \, \dd \mu \right | \ll  e^{-\kappa_{\mathrm{mix}} |t|} \; \cS^K_m(\phi_1) \, \cS^K_m(\phi_2).
$$
\end{fact}

\section{Reduction to a Haar deviation estimate} \label{sec:counting}
In this section we deduce Theorem~\ref{thm:Grassmann} from a Haar deviation bound (Proposition~\ref{prop:Haar}). We first relate the counting function to a Grassmann Siegel transform and use a Borel--Cantelli argument to reduce the desired almost-sure counting estimate to a deviation bound on $X_\ell$. We then deduce this deviation bound from the Haar deviation estimate in Proposition~\ref{prop:Haar} by a thickening argument.

\subsection{Counting lattice points} \label{sec:Volume}
For a non-increasing $\psi : [1,+\infty) \rightarrow (0,+\infty)$ and $T \geq 1$ we set
\begin{equation} \label{eq:Definition_E}
\cE_{\psi}(T) := \left \{ \bw \in \widetilde{X}_{k} : d([\bw], S_0) < \psi(\|\bw\|), \, 1 \leq \|\bw\| < T \right \} ,
\end{equation}
where, as before, $S_0 = \{z \in X_k : z \subseteq x_0 \}$. By the polar integration formula~\eqref{eq:polar_measure}, its volume agrees with the function introduced in Section~\ref{sec:main_result}:
\begin{equation}\label{eq:Volume_Exact}
\lambda(\cE_\psi(T))=2d\,a_{k,d}
\int_1^T v(\psi(t))t^{d-1}\,\dd t=M_\psi(T).
\end{equation}
Fix a Borel section $x \mapsto k_x$ of the map $K \rightarrow X_{\ell}$, $k \mapsto k^{-1}x_0$, so that $k_x^{-1} x_0 = x$.

\begin{lemma} \label{lem:counting_lattice_points}
For every $x \in X_{\ell}$ and $T \geq 1$,
\begin{align}\label{equ:coufor}
\cN_\psi(x,T) = \frac 12 \, \widehat{\mathbbm{1}}_{\cE_{\psi}(T)} (k_x \G).
\end{align}
\end{lemma}

\begin{proof}
By the $K$-invariance of $d(\cdot,\cdot)$ and $\|\cdot\|$ and since $k_x^{-1} S_0 = S(x)$, we have 
\begin{align*}
k_x^{-1} \cE_{\psi}(T) = \{ \bw\in \widetilde{X}_{k} : d([\bw], S(x)) < \psi(\|\bw\|), \, 1 \leq \|\bw\| < T \}.
\end{align*}
Using the identity $d(v,x)=d(v, S(x))$ in~\eqref{eq:distance_Schubert}, a rational point $v\in X_{k}(\Q)$ satisfies $d(v , x ) < \psi(H(v))$ and $H(v) < T$ if and only if both of its primitive representatives $\pm \bv \in \cP_{k}$ lie in this set, since $H(v)=\|\pm \bv\| \geq 1$.
We conclude that
\begin{displaymath}
\cN_\psi(x,T)=\frac12\sum_{\bv \in \cP_{k}}\mathbbm{1}_{k_x^{-1}\cE_{\psi}(T)}(\bv) =\frac{1}{2} \, \widehat{\mathbbm{1}}_{k_x^{-1}\cE_{\psi}(T)} (\G)=\frac{1}{2} \, \widehat{\mathbbm{1}}_{\cE_{\psi}(T)} (k_x \G). \qedhere
\end{displaymath}
\end{proof}

\subsection{Deviation bound}
We now state the deviation estimate on $X_{\ell}$, which implies Theorem~\ref{thm:Grassmann} and is deduced from the key Haar deviation bound in Proposition~\ref{prop:Haar}.

\begin{proposition} \label{prop:discrepancy}
In the setting of Theorem~\ref{thm:Grassmann}, there are constants $\gamma_0 = \gamma_0(k,\ell,d) \in (0,1)$ and $C_1 = C_1(k,\ell,d) > 0$, both independent of $\psi$, such that for every $\gamma \in (0, \gamma_0)$ there is $T_0 = T_0(k,\ell,d,\psi,\gamma) \geq 1$ with
\begin{equation} \label{eq:discrepancy}
\sigma \Big ( \big \{ x \in X_{\ell} : \big | \widehat{\mathbbm{1}}_{\cE_{\psi}(T)}(k_x \G) - M_{\psi}(T) \big | \geq C_1 \, M_{\psi}(T)^{1-\gamma} \big \} \Big ) \ll M_{\psi}(T)^{-\gamma_0}
\end{equation}
for all $T \geq T_0$.
\end{proposition}

\begin{proof}[Proof of Theorem~\ref{thm:Grassmann} assuming Proposition~\ref{prop:discrepancy}]
Fix an approximating function $\psi$ as in Theorem~\ref{thm:Grassmann}. Since $\psi>0$, $M_\psi$ is continuous and strictly increasing. Moreover, $M_\psi(T)\asymp_{k,\ell,d}\Psi(T)\to\infty$, so $M_\psi$ maps $[1,\infty)$ bijectively onto $[0,\infty)$.

Set $\gamma=\gamma_0/2$ and, for $n\geq1$, let
$$
m_n:=n^{1/\gamma},\qquad T_n:=M_\psi^{-1}(m_n)
\quad\text{and}\quad
Z_x(T):=\widehat{\mathbbm{1}}_{\cE_\psi(T)}(k_x\Gamma).
$$
By Proposition~\ref{prop:discrepancy}, we have for all sufficiently large $n$,
$$
\sigma\bigl(\{x:|Z_x(T_n)-m_n|\geq C_1m_n^{1-\gamma}\}\bigr)
\ll M_{\psi}(T_n)^{-\gamma_0}=m_n^{-\gamma_0}=n^{-2}.
$$
These measures are summable in $n$. Hence, by the Borel--Cantelli lemma, for $\sigma$-almost every $x$, we have
$|Z_x(T_n)-m_n|<C_1m_n^{1-\gamma}$ for all sufficiently large $n$.
Moreover, note that
$$
m_{n+1}-m_n\ll n^{1/\gamma-1}=m_n^{1-\gamma},
\quad\text{and}\quad m_{n+1}\asymp m_n.
$$
For such an $x$ and sufficiently large $T$, choose $n$ with $m_n\leq m:=M_\psi(T)<m_{n+1}$. Since $\mathcal{E}_{\psi}(T_1)\subseteq\mathcal{E}_{\psi}(T_2)$ for $T_1\leq T_2$, a standard interpolation argument gives
$$
\begin{aligned}
|Z_x(T)-m| &\leq\max\{|Z_x(T_n)-m_n|,|Z_x(T_{n+1})-m_{n+1}|\} +(m_{n+1}-m_n) \ll m^{1-\gamma}.
\end{aligned}
$$
By Lemma~\ref{lem:counting_lattice_points}, this bound implies, for $\sigma$-almost every $x$,
$$
\cN_\psi(x,T)=\tfrac12 Z_x(T) =\tfrac12M_\psi(T)+O_{x,\psi}\bigl(M_\psi(T)^{1-\gamma}\bigr),
$$
which proves~\eqref{eq:asymptotic}.
\end{proof}

\subsection{Transfer to Haar measure} \label{sec:transfer}
In this subsection, we bound the deviation in~\eqref{eq:discrepancy} with respect to $\sigma$ by deviations with respect to $\mu$. To do so, we thicken $K\Gamma$ in directions belonging to $P=\Stab_G(x_0)$. The counting regions $\mathcal E_\psi(T)$ are stable under small perturbations in these directions, in the sense that the perturbed regions lie between slightly enlarged and contracted ones. This allows us to control the original counting function from both above and below.

To obtain upper and lower bounds for the counting function, let $\mathcal F$ denote the class of all non-increasing functions $\psi:(0,\infty)\to(0,2r_*]$. This class contains the approximation functions in Theorem~\ref{thm:Grassmann} after extending them constantly to $(0,1)$. For $\psi\in\mathcal F$, $a_0\in[1/2,2]$ and $T>a_0$, let
\begin{align*}
\cE_\psi^{a_0}(T):=\{\bw\in\widetilde X_k:d([\bw],S_0)<\psi(\|\bw\|),\ a_0\leq\|\bw\|<T\}  
\end{align*}
and $M_\psi^{a_0}(T):=\lambda(\cE_\psi^{a_0}(T))$. By~\eqref{eq:polar_measure},
$$
M_\psi^{a_0}(T)=2d\,a_{k,d}\int_{a_0}^T v(\psi(t))t^{d-1}\,\dd t.
$$
Set $C_{\rm geom}:=2^{2k}-1$. For $0<\rho\leq\frac{1}{4C_{\rm geom}}$, set $\nu_\pm:=1\pm C_{\rm geom}\rho$ and
\begin{equation}\label{eq:psi_pm}
\psi_\pm(t):=\nu_\pm\psi(t/\nu_\pm),\qquad
\cE_\rho^+(T):=\cE_{\psi_+}^{\nu_-}(\nu_+T),\qquad
\cE_\rho^-(T):=\cE_{\psi_-}^{\nu_+}(\nu_-T).
\end{equation}
If $\psi\in\mathcal F$ and $\psi\leq r_*$, then $\psi_\pm\in\mathcal F$, since $\nu_\pm\in[3/4,5/4]$ and the functions $\psi_\pm$ are positive, non-increasing and bounded by $2r_*$.

The following proposition transfers the deviation estimate with respect to $\sigma$ to deviation estimates with respect to $\mu$ for the enlarged and contracted counting regions.  

\begin{proposition}\label{prop:transfer_Haar}
Let $\alpha_d:=\dim G-\dim K=\tfrac{d(d+1)}2-1$. There are constants $C_{\rm vol}\geq1$ and $\rho_1>0$ such that, for every $\psi\in\mathcal F$ with $\psi\leq r_*$, $T\geq3$ with $M_\psi(T)\geq1$, $0<\rho\leq\rho_1$ and $u\geq4C_{\rm vol}\rho M_\psi(T)$, we have 
\begin{align*}
\sigma\left(\left\{x:\left|\widehat{\mathbbm{1}}_{\cE_\psi(T)}(k_x\Gamma)-M_\psi(T)\right|\geq u\right\}\right)
&\ll\rho^{-\alpha_d}\sum_\pm
\mu\left(\left\{\left|\widehat{\mathbbm{1}}_{\cE_\rho^\pm(T)}-\lambda(\cE_\rho^\pm(T))\right|\geq u/2\right\}\right).
\end{align*}
\end{proposition}
To prove Proposition~\ref{prop:transfer_Haar}, we first construct a probability measure on a small slice in $P$ and compare its average over $K$ with Haar measure. 
Let $\pi : G\to G/K$ be the quotient map and $o=K\in G/K$. Let $m_{G/K}$ be the $G$-invariant measure on $G/K$ normalized so that, for every $f\in C_c(G)$,
\begin{equation}\label{eq:quotient_integration_GK}
\int_G f(g)\,\dd m_G(g)
=\int_{G/K}\int_K f(gk)\,\dd m_K(k)\,\dd m_{G/K}(gK).
\end{equation}
Here recall that $m_K$ is the Haar probability measure on $K$.
Since $G=PK$, the identification $P/(P\cap K)\simeq G/K$ gives a smooth local section $\mathfrak{s}: V \to P$ near $o$, with $\pi\circ \mathfrak{s}= \mathrm{id}$ and $\mathfrak{s}(o) = I_d$. Here $V$ is an open neighborhood of $o$. 
For $\rho>0$, define
$$
V_\rho:=\{gK\in G/K:\max\{\|g\|_{\mathrm{op}},\|g^{-1}\|_{\mathrm{op}}\}<1+\rho\} .
$$
This definition is independent of the representative $g$, since $K$ acts isometrically on $\R^d$. Choose $\rho_2\in(0,1]$ sufficiently small that $\overline{V}_{\rho_2}$ is compact and contained in $V$. For $0<\rho\leq\rho_2$, set
$$
S_\rho:=\mathfrak{s}(V_\rho),\qquad
m_\rho:=\mathfrak{s}_*\left(\frac{m_{G/K}|_{V_\rho}}{m_{G/K}(V_\rho)}\right).
$$
Here $m_\rho$ is the probability measure on $S_\rho$ obtained by pushing forward the normalized restriction of $m_{G/K}$ to $V_\rho$ under $\mathfrak{s}$.
\begin{lemma}\label{lem:thickening}
For $0<\rho\leq\rho_2$ and every non-negative Borel function $F$ on $\Omega$, we have
$$
\int_K\int_{S_\rho}F(pk\Gamma)\,\dd m_\rho(p)\,\dd m_K(k)
\ll\rho^{-\alpha_d}\int_\Omega F\,\dd\mu.
$$
\end{lemma}

\begin{proof}
Set $\mathfrak{C}:=\mathfrak{s}(\overline{V}_{\rho_2})K$, which is compact by the choice of $\rho_2$. For $0<\rho\leq\rho_2$, by~\eqref{eq:quotient_integration_GK} and the definition of $m_{\rho}$ we have
$$
\begin{aligned}
\int_K\int_{S_\rho}F(pk\Gamma)\,\dd m_\rho(p)\,\dd m_K(k)
&=\frac1{m_{G/K}(V_\rho)}\int_{V_\rho}\int_K
F(\mathfrak{s}(z)k\Gamma)\,\dd m_K(k)\,\dd m_{G/K}(z)\\
&=\frac1{m_{G/K}(V_\rho)}\int_{\pi^{-1}(V_\rho)}F(g\Gamma)\,\dd m_G(g).
\end{aligned}
$$
Since $\mathfrak{C}$ is compact, it can be covered by finitely many measurable fundamental domains for $G/\Gamma$, with their number $N$ independent of $\rho$ and $F$. Recall that $m_G$ agrees locally with $\mu$ under the quotient map $G\to\Omega$. Thus we have
$$
\int_{\mathfrak{C}}F(g\Gamma)\,\dd m_G(g)
\leq N\int_\Omega F\,\dd\mu.
$$
It thus remains to show that 
\begin{align}
m_{G/K}(V_{\rho})\asymp \rho^{\alpha_d},\qquad \forall\, 0<\rho\leq \rho_2.
\end{align}
For this, we invoke the polar decomposition in $G$. Let 
$$
\mathfrak p:=\{X\in \mathrm{M}_{d\times d}(\R):X^{\mathsf T}=X,\ \operatorname{tr}X=0\}.
$$
Then every $g\in G$ can be written uniquely in the form $g=e^Xk$, with $X\in\mathfrak p$ and $k\in K$.
The polar map $\Phi(X,k)=e^Xk$ is a smooth diffeomorphism, so Haar measure has a smooth positive density in these coordinates. Right $K$-invariance makes this density independent of $k$, giving
$$
\dd m_G(e^Xk)=J(X)\,\dd X\,\dd m_K(k),
$$
where $\dd X$ is Lebesgue measure on the $\alpha_d$-dimensional space $\mathfrak p$. By continuity, $J$ is bounded above and bounded away from zero near $0$.
Since $X$ is symmetric with real eigenvalues $\lambda_1,\ldots,\lambda_d$, we have $\|X\|_{\mathrm{op}}=\max_i|\lambda_i|$. Orthogonal multiplication preserves the operator norm, so
$$
\max\{\|e^Xk\|_{\mathrm{op}},\|(e^Xk)^{-1}\|_{\mathrm{op}}\}
=\max\{e^{\max_i\lambda_i},e^{-\min_i\lambda_i}\}
=e^{\|X\|_{\mathrm{op}}}.
$$
Thus, up to reducing $\rho_2$ if necessary, the quotient integration formula~\eqref{eq:quotient_integration_GK} and $m_K(K)=1$ give
\begin{align*}
m_{G/K}(V_\rho) =m_G(\pi^{-1}(V_\rho)) =\int_{\|X\|_{\mathrm{op}}<\log(1+\rho)}J(X)\,\dd X \asymp\bigl(\log(1+\rho)\bigr)^{\alpha_d}\asymp\rho^{\alpha_d}.
\end{align*}
This completes the proof.
\end{proof}

We next compare the regions under the action of the small parabolic slice $S_\rho$.

\begin{lemma}\label{lem:slice_comparison}
There is a constant $C_{\rm vol}\geq1$ such that the following holds for
$0<\rho\leq\rho_1:=\min\{\rho_2,(4C_{\rm geom})^{-1}\}$.
For $\psi\in\mathcal F$ with $\psi\leq r_*$ and $T\geq3$, let $\psi_\pm$ and $\cE_\rho^\pm(T)$ be as in~\eqref{eq:psi_pm}. Then, for every $p\in S_\rho$ and $x\in\Omega$, we have
\begin{align}\label{eq:slice_comparison1}
\widehat{\mathbbm{1}}_{\cE_\rho^-(T)}(px)
\leq\widehat{\mathbbm{1}}_{\cE_\psi(T)}(x)
\leq\widehat{\mathbbm{1}}_{\cE_\rho^+(T)}(px),
\end{align}
and
\begin{align}\label{eq:slice_comparison2}
|\lambda(\cE_\rho^\pm(T))-M_\psi(T)|\leq C_{\rm vol}\rho(M_\psi(T)+1).
\end{align}
\end{lemma}

\begin{proof}
Recall that $C_{\rm geom}=2^{2k}-1$; by convexity of the function $f(t)=(1+t)^{2k}$ on $[0,1]$, we have 
\begin{align}\label{equ:cgemobd}
(1+\rho)^{2k}\leq (1-\rho)f(0)+\rho f(1)=1+C_{\rm geom}\rho,\qquad \forall\, 0<\rho\leq 1.
\end{align}
Now fix $0<\rho\leq\rho_1$ and recall that $\nu_\pm=1\pm C_{\rm geom}\rho$. 
For $p\in S_\rho$, by the definition of $V_\rho$,~\eqref{eq:exterior_operator_norm},~\eqref{equ:cgemobd} and 
$\nu_-\leq\nu_+^{-1}$ we get
\begin{equation}\label{eq:slice_norm_comparison}
\nu_-\|\bw\|\leq\nu_+^{-1}\|\bw\|\leq\|p\bw\|\leq\nu_+\|\bw\|,\qquad \forall\, \bw\in\wedge^k\R^d.
\end{equation}
Since $p$ fixes $x_0$, it preserves $S_0$. Using $d(z,x_0)=d(z,S_0)$ from~\eqref{eq:distance_Schubert} and applying~\eqref{eq:grassmann_lipschitz_bound} to $p$ and $p^{-1}$, together with~\eqref{equ:cgemobd} and $\nu_-\leq\nu_+^{-1}$, we obtain
\begin{equation}\label{eq:proof_slice_comparison}
\nu_-d([\bw], S_0)\leq d([p\bw],S_0)\leq\nu_+d([\bw],S_0),\qquad \forall\, \bw\in\widetilde X_k.
\end{equation}
Thus if $\bw\in\cE_\psi(T)$, then $\nu_-\leq \|p\bw\|<\nu_+T$. Moreover, since $\psi$ is non-increasing, we have by~\eqref{eq:slice_norm_comparison} and~\eqref{eq:proof_slice_comparison}
$$
d([p\bw],S_0)\leq \nu_+d([\bw],S_0)<\nu_+\psi(\|\bw\|) \leq \nu_+\psi(\|p\bw\|/\nu_+)=\psi_+(\|p\bw\|).
$$
Hence $p\bw\in\cE_\rho^+(T)$. Conversely, if $p\bw\in\cE_\rho^-(T)$, then $\nu_+\leq \|p\bw\|<\nu_-T$, so~\eqref{eq:slice_norm_comparison} implies $1\leq \|\bw\|<T$. By~\eqref{eq:slice_norm_comparison}, \eqref{eq:proof_slice_comparison} and monotonicity of $\psi$, we have
$$
d([\bw],S_0)<\nu_-^{-1}\psi_-(\|p\bw\|)=\psi(\|p\bw\|/\nu_-)\leq\psi(\|\bw\|),
$$
so $\bw\in\cE_\psi(T)$. Taking Grassmann Siegel transforms of these indicator inequalities proves~\eqref{eq:slice_comparison1}.

To prove~\eqref{eq:slice_comparison2}, we apply the polar integration formula~\eqref{eq:polar_measure} and make the changes of variables $t=\nu_\pm s$ to obtain
\begin{align*}
\lambda(\cE_\rho^+(T))=2d\,a_{k,d}\nu_+^d
\int_{\nu_-/\nu_+}^T v(\nu_+\psi(s))s^{d-1}\,\dd s, \quad
\lambda(\cE_\rho^-(T))=2d\,a_{k,d}\nu_-^d
\int_{\nu_+/\nu_-}^T v(\nu_-\psi(s))s^{d-1}\,\dd s.
\end{align*}
By~\eqref{eq:RiemannianVolumeDelta}, we have
$$
v(\nu_\pm\psi(s))=v(\psi(s))(1+O(\rho)),
\qquad \nu_\pm^d=1+O(\rho).
$$
Here the implied constants depend only on $d$. The intervals between $1$ and $\nu_\mp/\nu_\pm$ have length $O(\rho)$ and lie in $[3/5,5/3]$, where $v\leq1$ and $s^{d-1}\ll1$. Therefore
$$
\lambda(\cE_\rho^\pm(T))=M_\psi(T)+O\bigl(\rho(M_\psi(T)+1)\bigr).
$$
This proves~\eqref{eq:slice_comparison2} by taking $C_{\rm vol}$ to be a bounding constant here. 
\end{proof}

\begin{proof}[Proof of Proposition~\ref{prop:transfer_Haar}]
Let $C_{\rm vol}$ and $\rho_1$ be as in Lemma~\ref{lem:slice_comparison}, and write $\cE=\cE_\psi(T)$ and $M=M_\psi(T)$. Recall that $X_\ell\simeq K/(K\cap P)$ via $k(K\cap P)\mapsto kx_0$. The map $k\mapsto k^{-1}x_0$ pushes $m_K$ forward to $\sigma$. Since $\cE$ is invariant under $K\cap P$, the value of $\widehat{\mathbbm{1}}_\cE(k\Gamma)$ depends only on $k^{-1}x_0\in X_{\ell}$, so
$$
\sigma\bigl(\{x:|\widehat{\mathbbm{1}}_\cE(k_x\Gamma)-M|\geq u\}\bigr)
=m_K\bigl(\{k:|\widehat{\mathbbm{1}}_\cE(k\Gamma)-M|\geq u\}\bigr).
$$
By~\eqref{eq:slice_comparison2}, since $M\geq1$ and $u\geq4C_{\rm vol}\rho M$, we have
$$
|\lambda(\cE_\rho^\pm(T))-M|\leq C_{\rm vol}\rho(M+1)\leq u/2.
$$
This, together with~\eqref{eq:slice_comparison1}, implies that for every $p\in S_\rho$, 
\begin{align*}
&\bigl\{k\in K:|\widehat{\mathbbm{1}}_\cE(k\Gamma)-M|\geq u\bigr\}
\subseteq\bigcup_\pm
\bigl\{k\in K:|\widehat{\mathbbm{1}}_{\cE_\rho^\pm(T)}(pk\Gamma)
-\lambda(\cE_\rho^\pm(T))|\geq u/2\bigr\}.
\end{align*}
Consequently, for every $k\in K$ and $p\in S_\rho$, we have
$$
\mathbbm{1}_{\{|\widehat{\mathbbm{1}}_\cE-M|\geq u\}}(k\Gamma)
\leq\sum_\pm
\mathbbm{1}_{\{|\widehat{\mathbbm{1}}_{\cE_\rho^\pm(T)}-\lambda(\cE_\rho^\pm(T))|\geq u/2\}}(pk\Gamma).
$$
Integrating the above inequality against the probability measures $m_\rho$ and $m_K$, and applying Lemma~\ref{lem:thickening}, proves the claim.
\end{proof}

\subsection{Deduction from the Haar deviation bound} \label{sec:haar_statement}\label{sec:proof_discrepancy}

In this subsection, we state the key Haar deviation bound and use it, together with Proposition~\ref{prop:transfer_Haar}, to prove Proposition~\ref{prop:discrepancy}.

\begin{proposition}\label{prop:Haar}
There are constants $B\geq1$ and $\kappa,c_0,\eta,\theta_0,C_0>0$, depending only on $k,\ell,d$, such that, for every $\psi\in\mathcal F$, $a_0\in[1/2,2]$ and $T\geq2a_0$ with $M:=\lambda(\cE_\psi^{a_0}(T))\geq2$, and for every $0<\theta\leq\theta_0$, $Q,R\geq1$ and $u\geq C_0\theta M$, we have
\begin{equation}\label{eq:Haar_deviation}
\mu\bigl(|\widehat{\mathbbm{1}}_{\cE_\psi^{a_0}(T)}-M|\geq u\bigr)
\ll\frac{\theta^{-B}RM}{u^2}
+\frac{\theta^{-1}M}{u}\bigl(R^{-\kappa}Q^{c_0}+Q^{-\eta}\bigr).
\end{equation}
The implied constant is uniform in these parameters.
\end{proposition}

\begin{proof}[Proof of Proposition~\ref{prop:discrepancy} assuming Proposition~\ref{prop:Haar}]
Take $B,\kappa,c_0,\eta,\theta_0,C_0$ from Proposition~\ref{prop:Haar}, and let $\alpha_d,C_{\rm vol},\rho_1$ be as in Proposition~\ref{prop:transfer_Haar}. Let 
\begin{align*}
C_1:=4\max\{1,C_{\rm vol},C_0\}
\quad\text{and}\quad
\rho_*:=\min\{\rho_1,(4C_{\rm vol})^{-1}\}.
\end{align*}
Fix $0<\gamma<\gamma_0$, where $\gamma_0\in(0,1)$ will be chosen below.
The approximation function in Theorem~\ref{thm:Grassmann} satisfies $\psi\in\mathcal F$ and $\psi\leq r_*$. Choose $T_0\geq4$ so that, for all $T\geq T_0$, we have $M:=M_\psi(T)\geq4$ and $M^{-\gamma}\leq\min\{\theta_0,\rho_*\}$. Set
$$
u:=C_1M^{1-\gamma}\quad\text{and}\quad \theta:=\rho:=M^{-\gamma},
$$
so that they satisfy $u\geq 4C_0\theta M$ and $u\geq 4C_{\rm vol}\rho M$. 
The comparison in Lemma~\ref{lem:slice_comparison} gives
$$
\tfrac12M\leq M_\pm:=\lambda(\cE_\rho^\pm(T))\leq2M.
$$
Hence we also have $u/2\geq C_0\theta M_\pm$. Applying Proposition~\ref{prop:Haar} to
$\cE_\rho^+(T)=\cE_{\psi_+}^{\nu_-}(\nu_+T)$ and
$\cE_\rho^-(T)=\cE_{\psi_-}^{\nu_+}(\nu_-T)$ with threshold $u/2$, and using Proposition~\ref{prop:transfer_Haar} and $M_\pm\asymp M$, we obtain, for any $R,Q\geq1$,
\begin{align*}
\sigma\bigl(\{x:|\widehat{\mathbbm{1}}_{\cE_\psi(T)}(k_x\Gamma)-M|\geq C_1M^{1-\gamma}\}\bigr)
&\ll \rho^{-\alpha_d}\left(\frac{\theta^{-B}RM}{u^2}
+\frac{\theta^{-1}M}{u}(R^{-\kappa}Q^{c_0}+Q^{-\eta})\right)\\
&\ll M^{(\alpha_d+B+2)\gamma}
\left(\frac RM+R^{-\kappa}Q^{c_0}+Q^{-\eta}\right).
\end{align*}
Now take $R=M^{1/2}$, $Q=M^{\kappa/(4c_0)}$ and
$$
\gamma_0:=\frac{1}{\alpha_d+B+3}
\min\left\{\frac12,\frac\kappa4,\frac{\kappa\eta}{4c_0}\right\}\in(0,1).
$$
Then the preceding bound is
\begin{align*}
\ll M^{(\alpha_d+B+2)\gamma} \left(M^{-1/2}+M^{-\kappa/4}+M^{-\kappa\eta/(4c_0)}\right) \ll M^{(\alpha_d+B+2)\gamma-(\alpha_d+B+3)\gamma_0} \leq M^{-\gamma_0},
\end{align*}
since $0<\gamma<\gamma_0$. As $\gamma_0$ depends only on $k,\ell,d$, this proves~\eqref{eq:discrepancy}.
\end{proof}

It remains to prove Proposition~\ref{prop:Haar}, which is the purpose of the remaining sections.

\section{Staircases, normalization and three estimates} \label{sec:staircase}
In this section, we approximate the region $\cE_\psi^{a_0}(T)$ by inner and outer staircases, uniformly in $\psi\in\mathcal F$ and $a_0\in[1/2,2]$. In particular, the construction applies to the regions $\cE_\rho^\pm(T)$ of the previous section. We normalize each shell with $a(1/r)$, which expands its transverse directions from scale $r$ to unit scale, and smoothly approximate the normalized shells. Finally, we decompose the resulting sums of Grassmann Siegel transforms as described in Section~\ref{sec:strategy} to reduce to the three estimates needed for the Haar deviation bound.

\subsection{Effective staircase approximation} \label{sec:mesh}

For $\psi\in\mathcal F$, $a_0\in[1/2,2]$ and $T>a_0$, let $\cE_\psi^{a_0}(T)$ and $M_\psi^{a_0}(T)$ be as in Section~\ref{sec:transfer}.
For $h,r>0$ and $\theta\geq0$, define the shell
$$
\cD(h,r,\theta) :=\left\{\bw\in\widetilde X_k: h\leq\|\bw\|<e^\theta h,\quad d([\bw],S_0)<r\right\}.
$$
The shell is empty when $\theta=0$. Recall that $v(r)$ denotes the $\sigma_k$-measure of the tube of radius $r$ around $S_0$. Polar integration, as in~\eqref{eq:Volume_Exact}, gives
\begin{equation}\label{eq:D_volume}
\lambda\bigl(\cD(h,r,\theta)\bigr) =2a_{k,d}\bigl(e^{d\theta}-1\bigr)h^d v(r).
\end{equation}

Assume that $T\geq2a_0$ and fix $0<\theta<\log 2$. Set
$$
h_j:=a_0e^{j\theta}\quad \text{and} \quad J:=\left\lfloor\theta^{-1}\log(T/a_0)\right\rfloor,
$$
so that $J\geq1$ and $h_J\leq T<h_{J+1}$. Since $\psi$ is non-increasing, on each height interval $[h_j,h_{j+1})$, the tube radius $\psi(\|\bw\|)$ lies between $\psi(h_{j+1})$ and $\psi(h_j)$. We therefore define the inner and outer staircases by
$$
\cE_\theta^- :=\bigcup_{j=0}^{J-1}\cD\bigl(h_j,\psi(h_{j+1}),\theta\bigr) \quad \text{and} \quad \cE_\theta^+ :=\bigcup_{j=0}^{J}\cD\bigl(h_j,\psi(h_j),\theta\bigr).
$$
The inner staircase stops at $h_J$, while the outer staircase extends to $h_{J+1}$. Thus
\begin{equation}\label{eq:staircase_inclusion}
\cE_\theta^-\subseteq\cE_\psi^{a_0}(T)\subseteq\cE_\theta^+.
\end{equation}

For simplicity, write $M=M_\psi^{a_0}(T)$.

\begin{lemma}\label{lem:mesh}
For every $T\geq2a_0$ with $M\geq2$, we have
$$
\lambda(\cE_\theta^\pm)=M+O(\theta M).
$$
\end{lemma}

\begin{proof}
For $0\leq j<J$, the inner shell at height $h_j$ and the outer shell at height $h_{j+1}$ have the same tube radius $\psi(h_{j+1})$. Since $h_{j+1}=e^\theta h_j$, formula~\eqref{eq:D_volume} shows that the outer shell has $e^{d\theta}$ times the volume of its inner counterpart. Summing over the pairs and adding the first outer shell gives
$$
\lambda(\cE_\theta^+) =\lambda\bigl(\cD(a_0,\psi(a_0),\theta)\bigr)+e^{d\theta}\lambda(\cE_\theta^-).
$$
Consequently,
\begin{align*}
\lambda(\cE_\theta^+)-\lambda(\cE_\theta^-) =(e^{d\theta}-1) \left(2a_{k,d}a_0^d v(\psi(a_0))+\lambda(\cE_\theta^-)\right) \ll \theta(1+M)\ll\theta M,
\end{align*}
where we used $a_0\leq2$, $v\leq1$, $M\geq2$, and $\lambda(\cE_\theta^-)\leq M$. Since $M$ lies between the two staircase volumes, the claim follows.
\end{proof}

\subsection{Shells and normalization} \label{sec:shells}

For $h>0$, $r>0$ and $\theta\geq0$, we define the normalized shell by
$$
\widetilde{\cD}(h,r,\theta):=a(1/r)\,\cD(h,r,\theta).
$$
For $|\tau|\leq1/2$, define a family of shells nested by inclusion as $\tau$ increases:
$$
\cD_\tau :=\Bigl\{ \bw \in \widetilde X_k : he^{-\tau}\leq\| \bw\|< he^{\theta+\tau}, d([\bw],S_0)<(1+\tau)r \Bigr\}, \qquad \widetilde{\cD}_\tau :=  a(1/r) \cD_\tau .
$$
Note that $\widetilde{\cD}_0 = \widetilde{\cD}(h,r,\theta)$. The following lemma describes the geometry of these normalized shells.

\begin{lemma} \label{lemma:shell_geometry}
Let $h > 0$, $ 0 < r \leq 2 r_*$, and $0 < \theta \leq 1$, and set $s := h r^{D/d}$. There exist $C \geq 1$ and $\delta_0 > 0$, depending on $k,\ell,d$, such that the following hold.
\begin{enumerate}[label=(\roman*)]
\item For every $|\tau| \leq 1/2$, 
$$
\widetilde{\cD}_\tau \subseteq \{\bw\in\widetilde X_k:s/C\leq\|\bw\|\leq Cs\} \quad \text{ and } \quad \lambda(\widetilde{\cD}_0)\asymp\theta s^d .
$$
\item For every $|\tau|\leq1/4$, $0<\delta \leq \delta_0$, and $g\in B_G(\delta)$, set $\eta:=C\delta$. Then
$$
\widetilde{\cD}_{\tau-\eta}\subseteq g\widetilde{\cD}_\tau \subseteq\widetilde{\cD}_{\tau+\eta}
\quad\text{and} \quad \lambda\bigl(\widetilde{\cD}_{\tau+\eta}\setminus \widetilde{\cD}_{\tau-\eta}\bigr) \ll \delta s^d.
$$
\end{enumerate}
\end{lemma} 
\begin{proof}
We first establish the required norm and distance estimates. Fix $0<\delta\leq\delta_0$, with $\delta_0$ to be chosen sufficiently small below, and let $g\in B_G(\delta)$. Set $g_r:=a(r)ga(1/r)$. We claim that uniformly for $\bw\in\widetilde X_k$ with $d([\bw],S_0)<2r$,
\begin{equation}\label{eq:normalized_shell_stability}
\|a(1/r)\bw\|\asymp r^{D/d}\|\bw\|,\qquad \left|\log\frac{\|g_r\bw\|}{\|\bw\|}\right|\ll\delta,\qquad \bigl|d([g_r\bw],S_0)-d([\bw],S_0)\bigr|\ll\delta r.
\end{equation}
Choose a sufficiently small fixed $0 < r_0 \leq 1$. For $r_0 \leq r \leq 2r_*$, the elements $a(r)$, $a(1/r)$ lie in a compact subset of $G$. Hence the first estimate in \eqref{eq:normalized_shell_stability} holds trivially, and $g_r\in B_G(C_1\delta)$ for some fixed $C_1 \geq 1$. By the standard local comparison of $d_G$ with the matrix operator norm near $I_d$, we have $\|g_r-I_d\|_{\mathrm{op}}\ll d_G(I_d,g_r)\ll\delta$. Hence $\|g_r\|_{\mathrm{op}}\leq1+O(\delta)$. By~\eqref{eq:exterior_operator_norm}, this implies $\|g_r\bw\|\leq\|\wedge^k g_r\|_{\mathrm{op}}\|\bw\|\leq e^{O(\delta)}\|\bw\|$; applying the same argument to $g_r^{-1}$ gives the corresponding lower bound. Combining these two bounds gives the second estimate in \eqref{eq:normalized_shell_stability}. For the third estimate, write $g_r=\exp Z$ with $\|Z\|\ll\delta$. For each fixed $z\in X_k$, smoothness of the action gives
$$
d(\exp(Z)z,z)\leq\int_0^1\left\|\frac{\dd}{\dd t}\exp(tZ)z\right\|\,\dd t\ll_z\|Z\|\ll \delta.
$$
The implied constant can be made uniform since $X_k$ is compact. Applying this estimate to $z=[\bw]$ and using triangle inequality, we get 
$$
|d([g_r\bw],S_0)-d([\bw],S_0)|\leq d(g_r[\bw],[\bw])\ll\delta\ll\delta r,
$$
since $r\geq r_0$. This proves the third estimate in \eqref{eq:normalized_shell_stability} for $r_0\leq r\leq 2r_*$.

Now for $0<r<r_0$, by~\eqref{eq:zooming_coordinates}, for $\bw\in\widetilde X_k$ with $d([\bw],S_0)<2r$ we may write
$$
\bw=c\,q\exp(ru)\be_0 \quad \text{with}\quad c\in\R^\times,\quad q\in H,\quad u\in \mathfrak{n} .
$$
Since $d([\bw],S_0)<2r$, the principal-angle identity~\eqref{eq:paidu}  bounds each singular value of $ru$ by $\tan(2r)$, so $\|ru\|_{\mathrm F}\leq\sqrt{k}\tan(2r)\leq1$ for $r_0$ sufficiently small. By~\eqref{eq:zooming_projective}, we have $\|r u\|_{\mathrm F} \asymp d(\exp(ru)z_0, S_0)=d([\bw],S_0)<2r$, implying that $\|u\|_{\mathrm F} < C_2$ for some $C_2 \ge 1$. Applying the cone zooming identity~\eqref{eq:zooming_cone} we get
$$
a(1/r)\bw=c\,r^{D/d}q\exp(u)\be_0.
$$
Set $J(u'):=\|\exp(u')\be_0\|$. In the expansion $\exp(u')\be_0=\bigwedge_{j=1}^k(e_j+u'e_j)$, the coefficient of $\be_0$ is $1$, so $J(u')\geq1$. This is a smooth bounded function on bounded subsets of $\mathfrak{n}$. Since $q$ acts isometrically, we have
$$
\|\bw\|=|c|J(ru)\asymp|c|\quad \text{and} \quad \|a(1/r)\bw\|=|c|r^{D/d}J(u)\asymp|c|r^{D/d}.
$$
This proves the first estimate.

We next perturb the zoomed vector $a(1/r)\bw$ by $g \in B_G(\delta)$ with $0 < \delta \leq \delta_0$. In fact, we claim that for $\delta_0$ sufficiently small, we may write
\begin{align}\label{eq:distorionide}
gq\exp(u)\be_0=\beta q_1\exp(u_1)\be_0\quad\text{with} \quad q_1\in H,\quad u_1\in \mathfrak{n},\quad \beta\in\R^\times,
\end{align}
with uniform bounds $|\beta-1|, \|u_1-u\|_{\mathrm F}\ll\delta$. Set $\widetilde g=q^{-1}gq$, so that $\|\widetilde g- I_d\|_{\mathrm{op}}\ll\delta$ uniformly in $q\in H$.
Let $\pi_{x_0}:\R^d\to\R^d$ be the orthogonal projection onto $x_0$. For $1\leq j\leq k$, we have $\pi_{x_0}(\exp(u)e_j)=e_j$, so $\|\pi_{x_0}(\widetilde g\exp(u)e_j)-e_j\|\ll\delta$ uniformly for $\|u\|_{\mathrm F}\leq C_2$. Since $e_1,\ldots,e_k$ form an orthonormal basis of $z_0$, the plane $z=\pi_{x_0}(\widetilde g\exp(u)z_0)$ spanned by these perturbed basis vectors satisfies $d(z,z_0)\ll\delta$. Choose a smooth local section $\mathfrak s$ of the orbit map $H\to S_0$, $h\mapsto hz_0$, near $z_0$, with $\mathfrak s(z_0)=I_d$ and $\mathfrak s(z)z_0=z$. As shown in Section~\ref{sec:zooming}, every element of $H$ commutes with $\pi_{x_0}$. Applying this to $\mathfrak s(z)^{-1}$, we obtain
$$
\pi_{x_0}(\mathfrak s(z)^{-1}\widetilde g\exp(u)z_0)
=\mathfrak s(z)^{-1}\pi_{x_0}(\widetilde g\exp(u)z_0)=z_0.
$$
By~\eqref{eq:graph_projection}, there are unique $u_1\in \mathfrak{n}$ and $\beta\in\R^\times$ such that
$$
\mathfrak s(z)^{-1}\widetilde g\exp(u)\be_0=\beta\exp(u_1)\be_0.
$$
This proves~\eqref{eq:distorionide} with $q_1=q\mathfrak s(z)$. Since $d(z,z_0)\ll\delta$ and $\mathfrak s$ is smooth, $\mathfrak s(z)=I_d+O(\delta)$ so that $\mathfrak s(z)^{-1}\widetilde g=I_d+O(\delta)$. 
This, together with the fact $\|u\|_{\mathrm F}\leq C_2$, implies that
$$
\|\beta\exp(u_1)\be_0-\exp(u)\be_0\|\ll\delta.
$$
In the same wedge expansion, each entry of $u:z_0\to x_0^\perp$ also appears, up to sign, as a coefficient in the standard wedge basis. Comparing the $\be_0$ coefficient and these entries gives $|\beta-1|,\|\beta u_1-u\|_{\mathrm F}\ll\delta$. Using also $\|u\|_{\mathrm{F}}\leq C_2$ we get $|\beta-1|\ll \delta$ and  $\|u_1-u\|_{\mathrm F}\ll\delta$ as desired.

Now, applying $a(r)$ to the identity in \eqref{eq:distorionide} and using the cone zooming identity \eqref{eq:zooming_cone} we get
$$
g_r\bw=c\beta q_1\exp(ru_1)\be_0.
$$
For $\delta_0$ sufficiently small, the bound on $\beta-1$ gives $\log\beta=O(\delta)$. Moreover, since $\|u\|_{\mathrm{F}}\leq C_2$, $\|u-u_1\|_{\mathrm{F}}\ll \delta$ and $0<r<r_0$, both $ru$ and $ru_1$ lie in a fixed ball in $\mathfrak{n}$. Then by the mean value theorem we have
$$
\left|\log\frac{\|g_r\bw\|}{\|\bw\|}\right| =\left|\log\beta+\log J(ru_1)-\log J(ru)\right| \ll\delta+r\|u_1-u\|_{\mathrm F} \ll\delta.
$$
Moreover, the map $u'\mapsto\exp(u')z_0$ has bounded differential on this fixed ball, so the image of the segment joining $ru$ to $ru_1$ has length $\ll r\|u_1-u\|_{\mathrm F}$. By the triangle inequality and using the facts that $q$ and $q_1$ preserve $S_0$ and act isometrically, we have
\begin{align*}
\bigl|d([g_r\bw],S_0)-d([\bw],S_0)\bigr| &=\bigl|d(\exp(ru_1)z_0,S_0)-d(\exp(ru)z_0,S_0)\bigr|\\ 
&\leq d(\exp(ru_1)z_0,\exp(ru)z_0)\\
&\ll r\|u_1-u\|_{\mathrm F} \ll\delta r.
\end{align*}
This completes the proof of~\eqref{eq:normalized_shell_stability}. We now use these estimates to prove (i) and (ii) of this lemma.

For $|\tau|\leq1/2$, every $\bw\in\cD_\tau$ satisfies $\|\bw\|\asymp h$ and $d([\bw],S_0)<3r/2$. Then by the first estimate in~\eqref{eq:normalized_shell_stability} we have
$$
\|a(1/r)\bw\|\asymp hr^{D/d}=s,
$$
proving the first assertion in (i) after choosing $C$ sufficiently large. For the second estimate in (i), by invariance of $\lambda$, the shell-volume formula~\eqref{eq:D_volume} and the $v(r)$ formula~\eqref{equ:vrcomp},
$$
\lambda(\widetilde{\cD}_0) =2a_{k,d}(e^{d\theta}-1)h^dv(r) \asymp\theta h^dr^D =\theta s^d.
$$

Next, for (ii), enlarge $C$ to dominate the constants in the last two estimates of~\eqref{eq:normalized_shell_stability}, and decrease $\delta_0$ so that $C\delta_0\leq1/4$. Fix $|\tau|\leq1/4$, let $\eta=C\delta$, and let $\bw\in\cD_\tau$. Then
$$
he^{-\tau-\eta}\leq\|g_r\bw\| <he^{\theta+\tau+\eta}\quad \text{and} \quad d([g_r\bw],S_0)<(1+\tau+\eta)r.
$$
Thus $g_r\cD_\tau\subseteq\cD_{\tau+\eta}$. Applying the same estimates to $g^{-1}$ and points of $\cD_{\tau-\eta}$ gives $g_r^{-1}\cD_{\tau-\eta}\subseteq\cD_\tau$.
Hence
$$
\cD_{\tau-\eta}\subseteq g_r\cD_\tau \subseteq\cD_{\tau+\eta}.
$$
Applying $a(1/r)$ proves the required inclusions for the normalized shells in (ii).

Finally, by invariance of $\lambda$ and the shell-volume formula~\eqref{eq:D_volume}, we have for $|\tau|\leq1/2$,
$$
\lambda(\widetilde{\cD}_\tau) =2a_{k,d}h^d \bigl(e^{d(\theta+\tau)}-e^{-d\tau}\bigr)_+ v((1+\tau)r),
$$
where $x_+:=\max\{x,0\}$. Since $|\tau|\leq1/4$ and $\eta/(1+\tau)\leq1/3$, by~\eqref{eq:RiemannianVolumeDelta} and~\eqref{equ:vrcomp} we have
$$
v((1+\tau\pm\eta)r)=v((1+\tau)r)(1+O(\eta))\ll r^D.
$$
Let $F(t):=(e^{d(\theta+t)}-e^{-dt})_+$. Since the shells are nested and $F$ is bounded on $[-1/2,1/2]$, we have
\begin{align*}
\lambda\bigl(\widetilde{\cD}_{\tau+\eta}\setminus\widetilde{\cD}_{\tau-\eta}\bigr)
&=2a_{k,d}h^dv((1+\tau)r)\bigl(F(\tau+\eta)-F(\tau-\eta)+O(\eta)\bigr)\\
&\ll h^dr^D\bigl(F(\tau+\eta)-F(\tau-\eta)+\eta\bigr).
\end{align*}
Since $s^d=h^dr^D$ and $\eta=C\delta$, it suffices to show $F(\tau+\eta)-F(\tau-\eta)\ll\eta$. This follows from the mean value theorem and $|x_+-y_+|\leq|x-y|$.
This proves the lemma.
\end{proof}

\subsection{Smoothing and regular targets}
\label{sec:staircase_conventions}

In this subsection, we construct smooth upper and lower approximations to the normalized shells by convolution, using Lemma~\ref{lemma:shell_geometry}.

The Lie algebra $\kg=\mathfrak{sl}_d(\R)$ acts on smooth functions on $\widetilde X_k$ by
$$
(Zf)(\bw) :=\left.\frac{\dd}{\dd t}f(\exp(tZ)\bw)\right|_{t=0}.
$$
As in Section~\ref{sec:Equidistribution}, we extend this action to monomials $\mathscr D$ in the chosen basis of $\kg$. Observe that for $f\in C_c^\infty(\widetilde X_k)$, we have $\mathscr D\widehat f=\widehat{\mathscr Df}$. 

Fix nonnegative smooth kernels $\rho_\vartheta$ on $G$, for $0<\vartheta\leq1$, satisfying
\begin{align}\label{eq:kernefun}
\supp\rho_\vartheta\subset B_G(\vartheta), \qquad \int_G\rho_\vartheta\,\dd m_G=1, \qquad \|\mathscr D\rho_\vartheta\|_{L^1(m_G)} \ll_{\mathscr D}\vartheta^{-\deg\mathscr D}.
\end{align}
For a bounded Borel function $F$ on $\widetilde X_k$, define the convolution 
$$
(\rho_\vartheta*F)(\bw) :=\int_G\rho_\vartheta(g)F(g^{-1}\bw)\,\dd m_G(g).
$$
By the convolution identity $\mathscr D(\rho_\vartheta*F)=(\mathscr D\rho_\vartheta)*F$ and the derivative bounds for $\rho_\vartheta$, we have
\begin{equation}\label{eq:convolution_derivative_bound}
\|\mathscr D(\rho_\vartheta*F)\|_\infty \leq\|\mathscr D\rho_\vartheta\|_{L^1(m_G)}\|F\|_\infty \ll_{\mathscr D} \vartheta^{-\deg\mathscr D}\|F\|_\infty.
\end{equation}

Let $C,\delta_0$ be the constants in Lemma~\ref{lemma:shell_geometry}, with $C\delta_0\leq1/4$. Fix a sufficiently large constant $C_*>2$, depending on $k,\ell,d$, as specified in the proof below.

\begin{definition}\label{def:regular}
Let $s>0$ and $0<\vartheta\leq1$. A smooth function $f:\widetilde X_k\to[0,1]$ is \emph{$(s,\vartheta)$-regular} if
\begin{enumerate}[label=(\roman*)]
\item $\supp f\subseteq \{\bw\in\widetilde X_k:s/C_*\leq\|\bw\|\leq C_*s\}$;
\item $\|\mathscr Df\|_\infty \ll_{\mathscr D}\vartheta^{-\deg\mathscr D}$ for every monomial $\mathscr D$;
\item $\|f(g^{-1}\cdot)-f\|_{L^1(\lambda)} \leq C_*\delta s^d$ for every $g\in B_G(\delta)$ and $0<\delta\leq\delta_0$.
\end{enumerate}
\end{definition}

\begin{lemma}\label{lem:shells_regular}
For $C_*$ sufficiently large, there exists $\vartheta_0\in(0,1]$ such that the following holds. Let $h>0$, $0<r\leq2r_*$, $0<\theta\leq1/10$, and let
$\widetilde{\cD}:=\widetilde{\cD}(h,r,\theta)$ and $s:=hr^{D/d}$.
For every $0<\vartheta\leq\vartheta_0$, there exist $(s,\vartheta)$-regular functions $f^\pm$ satisfying
\begin{align}\label{eq:smfunprop}
0\leq f^-\leq\mathbbm{1}_{\widetilde{\cD}} \leq f^+\leq1
\quad\text{and} \quad \lambda(f^+-f^-)\ll\vartheta s^d.
\end{align}
\end{lemma}

\begin{proof}
Choose $\vartheta_0\leq\min\{1,\delta_0\}$ and let $\eta=C\vartheta$, so that $\eta\leq1/4$. Using the family $\widetilde{\cD}_\tau$ from Lemma~\ref{lemma:shell_geometry}, define the smooth functions 
$$
f^\pm:=\rho_\vartheta* \mathbbm{1}_{\widetilde{\cD}_{\pm\eta}}.
$$
Then clearly we have $0 \leq f^\pm \leq 1$.
For $g\in B_G(\vartheta)$, part (ii) of Lemma~\ref{lemma:shell_geometry}, applied at $\tau=0$ to $g^{-1}$, gives
$$
\widetilde{\cD}_{-\eta} \subseteq g^{-1}\widetilde{\cD}_0 \subseteq\widetilde{\cD}_{+\eta}.
$$
Consequently,
$$
\mathbbm{1}_{\widetilde{\cD}_{-\eta}}(g^{-1}\bw) \leq\mathbbm{1}_{\widetilde{\cD}_0}(\bw) \leq\mathbbm{1}_{\widetilde{\cD}_{+\eta}}(g^{-1}\bw).
$$
Integrating against $\rho_\vartheta$ proves the pointwise bounds in~\eqref{eq:smfunprop}. By invariance of $\lambda$ and the second estimate in Lemma~ \ref{lemma:shell_geometry}~(ii),
$$
\lambda(f^+-f^-) =\lambda\bigl( \widetilde{\cD}_{+\eta}\setminus \widetilde{\cD}_{-\eta}\bigr) \ll\vartheta s^d.
$$

Next, we verify $f^{\pm}$ are both $(s,\vartheta)$-regular. By part (i) of Lemma~\ref{lemma:shell_geometry}, both $\widetilde{\cD}_{\pm\eta}$ lie in a fixed annulus of scale $s$. Elements of $B_G(\vartheta_0)$ change norms in $\wedge^k\R^d$ by uniformly bounded factors. Thus the supports of $f^\pm$ lie in $\{s/C_*\leq\|\bw\|\leq C_*s\}$ after choosing $C_*$ sufficiently large. The derivative bounds follow from~\eqref{eq:convolution_derivative_bound}.

For the translation bound, observe that by part (ii) of Lemma~\ref{lemma:shell_geometry} and the nesting of the shell family,
$$
\lambda\bigl(b\widetilde{\cD}_\tau \triangle\widetilde{\cD}_\tau\bigr) \leq \lambda\bigl(\widetilde{\cD}_{\tau+C\varepsilon}\setminus \widetilde{\cD}_{\tau-C\varepsilon}\bigr)\ll\varepsilon s^d
$$
whenever $|\tau|\leq1/4$, $b\in B_G(\varepsilon)$, and $0<\varepsilon\leq\delta_0$. Also, the volume and boundary estimates at $\tau=0$ give
$$
\lambda(\widetilde{\cD}_{\pm\eta}) \ll(\theta+\vartheta)s^d\ll s^d.
$$

Let $g\in B_G(\delta)$ with $0<\delta\leq\delta_0$. By invariance of $\lambda$ and Fubini's theorem,
\begin{align*}
\|f^\pm(g^{-1}\cdot)-f^\pm\|_{L^1(\lambda)} &\leq\int_G\rho_\vartheta(q)\, \lambda\bigl( (q^{-1}gq)\widetilde{\cD}_{\pm\eta} \triangle\widetilde{\cD}_{\pm\eta} \bigr)\,\dd m_G(q).
\end{align*}
Compact conjugation gives a constant $C_{\mathrm{conj}}\geq1$ such that $q^{-1}gq\in B_G(C_{\mathrm{conj}}\delta)$ uniformly for $q\in\supp\rho_\vartheta\subset B_G(1)$. If $C_{\mathrm{conj}}\delta\leq\delta_0$, the preceding translation estimate, applied at $\tau=\pm\eta$, bounds the integrand by $O(\delta s^d)$. If $C_{\mathrm{conj}}\delta>\delta_0$, the same bound follows from
$$
\lambda\bigl( (q^{-1}gq)\widetilde{\cD}_{\pm\eta} \triangle\widetilde{\cD}_{\pm\eta} \bigr) \leq2\lambda(\widetilde{\cD}_{\pm\eta}) \ll s^d\ll\delta s^d.
$$
Since $\rho_\vartheta$ has integral one, we conclude that
$$
\|f^\pm(g^{-1}\cdot)-f^\pm\|_{L^1(\lambda)} \ll\delta s^d.
$$
Enlarging $C_*$ if necessary, we obtain (iii), which completes the proof.
\end{proof}

\subsection{Smooth counts}\label{sec:smooth_counts}
We apply the smoothing lemma to the staircase shells.
Let
$$
r_j^+:=\psi(h_j)\quad(0\leq j\leq J)\quad \text{and} \quad r_j^-:=\psi(h_{j+1})\quad(0\leq j<J)
$$
be the radii of the upper and lower staircase shells, respectively. Since $\psi\in\mathcal F$, we have $r_j^\pm\leq2r_*$. Define the normalization times by
\begin{equation}\label{eq:times}
t_j^\pm:=-\log r_j^\pm.
\end{equation}
Each sequence $(t_j^\pm)$ is non-decreasing because $\psi$ is non-increasing. By Lemma~\ref{lemma:shell_geometry} (i), the vectors in $\widetilde{\cD}(h_j,r_j^\pm,\theta)$ have norms comparable to
\begin{equation}\label{eq:times_2}
s_j^\pm :=h_j(r_j^\pm)^{D/d} =a_0e^{j\theta-(D/d)t_j^\pm}.
\end{equation}
We call $s_j^\pm$ the normalized height of the shell. These heights may be arbitrarily small or large. When one staircase is fixed, we suppress the superscript $\pm$. By the definition of the Grassmann Siegel transform and \eqref{eq:staircase_inclusion}, we have for every $x \in \Omega$,
\begin{equation}\label{eq:upper_lower_staircase_siegel_transform}
\sum_{j=0}^{J-1} \widehat{\mathbbm{1}}_{\widetilde{\cD}(h_j,r_j^-,\theta)} \bigl(a(e^{t_j^-})x\bigr) \leq \widehat{\mathbbm{1}}_{\cE_\psi^{a_0}(T)}(x) \leq \sum_{j=0}^{J} \widehat{\mathbbm{1}}_{\widetilde{\cD}(h_j,r_j^+,\theta)}
\bigl(a(e^{t_j^+})x\bigr).
\end{equation}
To apply the dynamical estimates in Section~\ref{sec:Equidistribution}, we smooth the normalized shell indicators at scale $\vartheta=\theta^2$. The following proposition shows that the resulting sums still bound the count from below and above, and that each has mean $M+O(\theta M)$.

\begin{proposition}\label{prop:staircase}
Set $\theta_0:=\min\{1/10,\sqrt{\vartheta_0}\}$. Let $\psi\in\mathcal F$, $a_0\in[1/2,2]$, $T\geq2a_0$, and $0<\theta\leq\theta_0$, and suppose $M:=M_\psi^{a_0}(T)\geq2$. Let $I^+:=\{0,\ldots,J\}$ and $ I^-:=\{0,\ldots,J-1\}$. For the times and heights defined in~\eqref{eq:times} and~\eqref{eq:times_2}, there exist functions $(f_j^\pm)_{j\in I^\pm}$ with the following properties.
\begin{enumerate}[label=(\roman*)]
\item
Each $f_j^\pm$ is $(s_j^\pm,\theta^2)$-regular.
\item
Let $\cS^\pm:=\sum_{j\in I^\pm}\widehat{f_j^\pm}\circ a(e^{t_j^\pm})$ denote the smooth counts. Then 
$
\cS^-\leq \widehat{\mathbbm{1}}_{\cE_\psi^{a_0}(T)} \leq\cS^+.
$
\item $\displaystyle \sum_{j\in I^\pm}(s_j^\pm)^d\ll\theta^{-1}M$.
\item $\displaystyle \int_\Omega\cS^\pm\,\dd\mu =\sum_{j\in I^\pm}\lambda(f_j^\pm) =M+O(\theta M)$.
\end{enumerate}
\end{proposition}

\begin{proof}
For each $j\in I^+$, let $f_j^+$ be the upper approximation to $\widetilde{\cD}(h_j,r_j^+,\theta)$ supplied by Lemma~\ref{lem:shells_regular} at scale $\vartheta=\theta^2$. For each $j\in I^-$, let $f_j^-$ be the corresponding lower approximation to $\widetilde{\cD}(h_j,r_j^-,\theta)$. This proves (i) and gives
$$
f_j^-\leq \mathbbm{1}_{\widetilde{\cD}(h_j,r_j^-,\theta)}
\quad\text{and} \quad \mathbbm{1}_{\widetilde{\cD}(h_j,r_j^+,\theta)} \leq f_j^+.
$$
Positivity of the Grassmann Siegel transform and~\eqref{eq:upper_lower_staircase_siegel_transform} imply (ii).

For (iii), Lemma~\ref{lemma:shell_geometry} (i) and invariance of $\lambda$ give
$$
\lambda\bigl(\cD(h_j,r_j^\pm,\theta)\bigr) \asymp\theta(s_j^\pm)^d.
$$
The shells of each staircase are disjoint, so Lemma~\ref{lem:mesh} yields
$$
\theta\sum_{j\in I^\pm}(s_j^\pm)^d \ll \sum_{j\in I^\pm} \lambda\bigl(\cD(h_j,r_j^\pm,\theta)\bigr) =\lambda(\cE_\theta^\pm) \ll M.
$$

Finally, for (iv), the pointwise bounds and mass estimate in Lemma~\ref{lem:shells_regular} imply
$$
\left| \lambda(f_j^\pm) -\lambda\bigl(\widetilde{\cD}(h_j,r_j^\pm,\theta)\bigr) \right| \ll\theta^2(s_j^\pm)^d.
$$
By (iii), the sum of these errors is $O(\theta M)$. The mean value formula~\eqref{eq:WeilIntegrationFormula*}, invariance of $\mu$ and $\lambda$, and Lemma~\ref{lem:mesh} therefore give
\begin{displaymath}
\int_\Omega\cS^\pm\,\dd\mu =\sum_{j\in I^\pm}\lambda(f_j^\pm) =\lambda(\cE_\theta^\pm)+O(\theta M) =M+O(\theta M).\qedhere
\end{displaymath}
\end{proof}

\subsection{Truncation and the three estimates}
\label{sec:three_estimates}
Finally, we use the smooth counts $\cS^\pm$ of Proposition~\ref{prop:staircase} to reduce the Haar deviation bound to three estimates. For shells of small normalized height, we truncate their Grassmann Siegel transforms and estimate the resulting variance and truncation error. For shells of large normalized height, we estimate the total discrepancy directly in $L^1(\mu)$.

To truncate the Grassmann Siegel transforms of the shell functions, fix $b:=\binom{d}{k}+1$. For $c>0$, define
\begin{equation}\label{equ:Wcdef}
W_c:=\widehat{\varphi_c}\quad \text{with} \quad \varphi_c(\bw):=(1+\|\bw\|^2/c^2)^{-b} \quad(\bw\in\widetilde X_k).
\end{equation}
The function $W_c$ is a smooth $K$-invariant weighted count of primitive pure vectors in $\wedge^k\R^d$. Its weight satisfies $\varphi_c(\bw)\geq2^{-b}$ for $\|\bw\|\leq c$ and tends to zero as $\|\bw\|/c\to\infty$.

If $f$ is $(c,\vartheta)$-regular, then $0\leq f\leq1$ and $\supp f\subseteq\{\|\bw\|\leq C_*c\}$, so
\begin{align}\label{eq:sieweicomp}
f\leq(1+C_*^2)^b\varphi_c \quad \text{and}\quad \widehat f\leq(1+C_*^2)^bW_c.
\end{align}
Thus the tail of $W_c$ controls the truncation error (see~\eqref{eq:truncation_error_pointwise}), while its $K$-invariance simplifies the Sobolev norm estimates used for the variance bound (see the proof of Lemma~\ref{lem:truncation}).

Fix $\chi\in C^\infty([0,\infty))$ with $0\leq\chi\leq1$, $\chi=1$ on $[0,1]$, and $\chi=0$ on $[2,\infty)$. For each smooth count $\cS^\pm$ of Proposition~\ref{prop:staircase}, write
$$
M_\pm:=\mu(\cS^\pm).
$$
Given $Q,R\geq1$, we split the shell indices according to their normalized heights:
$$
I_{\mathrm{small}}^\pm :=\{j\in I^\pm:s_j^\pm\leq Q\}
\quad\text{and}\quad 
I_{\mathrm{large}}^\pm :=\{j\in I^\pm:s_j^\pm>Q\}.
$$
For $j\in I_{\mathrm{small}}^\pm$, define the truncated Grassmann Siegel transform
$$
\Phi_j^\pm :=\widehat{f_j^\pm}\, \chi(W_{s_j^\pm}/R).
$$
We then have the decomposition
\begin{equation}\label{eq:decomposition}
\cS^\pm-M_\pm = \bigl(\cS_0^\pm-\mu(\cS_0^\pm)\bigr) +\bigl(\cA^\pm-\mu(\cA^\pm)\bigr) +\cB^\pm,
\end{equation}
where
\begin{align*}
\cS_0^\pm =\sum_{j\in I_{\mathrm{small}}^\pm} \Phi_j^\pm\circ a(e^{t_j^\pm}), \,\,\,\, \cA^\pm =\sum_{j\in I_{\mathrm{small}}^\pm} \bigl(\widehat{f_j^\pm}-\Phi_j^\pm\bigr) \circ a(e^{t_j^\pm}), \,\,\,\, \cB^\pm =\sum_{j\in I_{\mathrm{large}}^\pm} \left( \widehat{f_j^\pm}\circ a(e^{t_j^\pm}) -\lambda(f_j^\pm) \right).
\end{align*}
Here $\cS_0^\pm$ is the truncated small-shell count, $\cA^\pm$ is its nonnegative truncation error, and $\cB^\pm$ is the centered large-shell count.

If $|\cS^\pm-M_\pm|\geq u/2$, then at least one of the three terms in~\eqref{eq:decomposition} has absolute value at least $u/6$. Applying Chebyshev's inequality to the first term and Markov's inequality to the other two, and using $\|\cA^\pm-\mu(\cA^\pm)\|_{1} \leq2\|\cA^\pm\|_{1}$, we obtain
\begin{equation}\label{eq:smooth_count_probability}
\mu\bigl(|\cS^\pm-M_\pm|\geq u/2\bigr) \leq \frac{36\Var_\mu(\cS_0^\pm)}{u^2} +\frac{12\|\cA^\pm\|_1+6\|\cB^\pm\|_1}{u}.
\end{equation}

We now reduce the proof of the Haar deviation bound in Proposition~\ref{prop:Haar} to the following three estimates.

\begin{proposition}\label{prop:three_bounds}
There exist constants $B\geq1$ and $\kappa,c_0,\eta>0$, depending only on $k,\ell,d$, such that, in the setting of Proposition~\ref{prop:staircase}, for every $Q,R\geq1$ and each sign,
\begin{enumerate}[label=(\roman*)]
\item $\|\cA^\pm\|_1 \ll R^{-\kappa}Q^{c_0}\theta^{-1}M$;
\item $\Var_\mu(\cS_0^\pm)
\ll\theta^{-B}RM$;
\item $\|\cB^\pm\|_1 \ll Q^{-\eta}\theta^{-1}M$.
\end{enumerate}
\end{proposition}
\begin{proof}[Proof of Proposition~\ref{prop:Haar}
assuming Proposition~\ref{prop:three_bounds}] 
By Proposition~\ref{prop:staircase} (ii) and (iv),
$$
\cS^- \leq\widehat{\mathbbm{1}}_{\cE_\psi^{a_0}(T)} \leq\cS^+ \quad \text{and}\quad |M_\pm-M|\leq C\theta M,
$$
where $C$ depends only on $k,\ell,d$. Choose $C_0\geq2C$ in Proposition~\ref{prop:Haar}. For $u\geq C_0\theta M$, we have $|M_\pm-M|\leq u/2$. Consequently,
\begin{align*}
\widehat{\mathbbm{1}}_{\cE_\psi^{a_0}(T)}\geq M+u \ \,\,\Longrightarrow \,\,
\cS^+-M_+\geq u/2,\qquad \text{and} \qquad \widehat{\mathbbm{1}}_{\cE_\psi^{a_0}(T)}\leq M-u
\,\, \Longrightarrow \,\, M_--\cS^-\geq u/2.
\end{align*}
Applying~\eqref{eq:smooth_count_probability} to each sign and using Proposition~\ref{prop:three_bounds} gives
\begin{align*}
\mu\left( \left|\widehat{\mathbbm{1}}_{\cE_\psi^{a_0}(T)}-M\right| \geq u \right) \leq \sum_\pm\mu\bigl(|\cS^\pm-M_\pm|\geq u/2\bigr) \ll \frac{\theta^{-B}RM}{u^2} +\frac{\theta^{-1}M}{u} \bigl(R^{-\kappa}Q^{c_0}+Q^{-\eta}\bigr).
\end{align*}
This proves~\eqref{eq:Haar_deviation}.
\end{proof}

\section{Moment and tail estimates for the Grassmann Siegel transform} \label{sec:tail}

In this section, we prove a moment bound for $W_c$ uniformly in the height $c$. The resulting tail estimate controls the truncation error in Proposition~\ref{prop:three_bounds}~(i).

We fix $\eps_0 \in (0,1]$ such that Fact~\ref{fact:L1}~(iv) holds with $\eps = \eps_0$, and set $\eps_1 := \eps_0/2$ and $p_1 := 1 + \eps_1$. For $t>0$, write $N_t:=\widehat{\mathbbm{1}}_{B_t}$.

\begin{proposition}\label{prop:tail}
For every $c>0$,
\begin{align}\label{eq:weighted_moment}
\int_\Omega W_c^{p_1}\,\dd\mu\ll c^d\max\{1,c\}^{d\eps_1}.
\end{align}
Consequently, for $R\geq1$,
\begin{align}\label{eq:tail}
\int_\Omega W_c\mathbbm{1}_{\{W_c>R\}}\,\dd\mu
\ll R^{-\eps_1}c^d\max\{1,c\}^{d\eps_1}.
\end{align}
\end{proposition}

\begin{remark}\label{rmk:weighted_moment_order}
The upper bound in~\eqref{eq:weighted_moment} is optimal. In fact, for every $c>0$, we have
$$
\int_\Omega W_c^{p_1}\,\dd\mu\asymp c^d\max\{1,c\}^{d\eps_1}.
$$
See Remark~\ref{rmk:weighted_moment_lower} below for more details.
\end{remark}

\subsection{Proof of Proposition~\ref{prop:three_bounds}~(i)} \label{sec:proof_bias}
Assuming Proposition~\ref{prop:tail}, we deduce Proposition~\ref{prop:three_bounds}~(i) as follows.

\begin{proof}
Fix either sign in Proposition~\ref{prop:staircase}, and write $\{f_j\}_{j\in I}$ for the corresponding smooth shell functions, with normalized heights $s_j$ and times $t_j$. Each $f_j$ is $(s_j,\theta^2)$-regular by Proposition~\ref{prop:staircase}~(i), so the comparison in~\eqref{eq:sieweicomp} gives $\widehat f_j\ll W_{s_j}$ uniformly in $j$.
Recall that $\Phi_j=\widehat{f_j}\chi(W_{s_j}/R)$, where $0\leq\chi\leq1$ and $\chi=1$ on $[0,1]$. Thus
\begin{align}\label{eq:truncation_error_pointwise}
0\leq\widehat{f_j}-\Phi_j
&=\widehat{f_j}\bigl(1-\chi(W_{s_j}/R)\bigr)
\leq\widehat{f_j}\mathbbm{1}_{\{W_{s_j}>R\}}
\ll W_{s_j}\mathbbm{1}_{\{W_{s_j}>R\}}.
\end{align}
By invariance of $\mu$ and the triangle inequality we have 
\begin{align*}
\|\cA\|_1
&\leq\sum_{j\in I_{\mathrm{small}}}
\bigl\|(\widehat{f_j}-\Phi_j)\circ a(e^{t_j})\bigr\|_1
=\sum_{j\in I_{\mathrm{small}}}\|\widehat{f_j}-\Phi_j\|_1\\
&\ll R^{-\eps_1}Q^{d\eps_1}
\sum_{j\in I_{\mathrm{small}}}s_j^d
\ll R^{-\eps_1}Q^{d\eps_1}\theta^{-1}M.
\end{align*}
Here the second last inequality follows from~\eqref{eq:truncation_error_pointwise} and Proposition~\ref{prop:tail}, since $0<s_j\leq Q$ for $j\in I_{\mathrm{small}}$, and the last inequality follows from Proposition~\ref{prop:staircase}~(iii).
Taking $\kappa=\eps_1$ and $c_0=d\eps_1$ proves Proposition~\ref{prop:three_bounds}~(i).
\end{proof}

Note that the tail estimate~\eqref{eq:tail} follows from~\eqref{eq:weighted_moment} since
$W_c\mathbbm{1}_{\{W_c>R\}}\leq R^{-\eps_1}W_c^{p_1}$. Thus it remains to prove the moment estimate~\eqref{eq:weighted_moment}.

\subsection{Annular counts and weighted counts in $L^{p_1}$} \label{sec:Lp_bounds}

To estimate the contribution to $W_c$ from vectors of norm greater than one, we first bound the count in an annulus.

\begin{lemma}\label{lem:annular_count}
For any $t\geq1$, let $A_t:=\{\bw\in\widetilde X_k:t<\|\bw\|\leq 2t\}$. Then
\begin{align}\label{eq:annular_count}
\|\widehat{\mathbbm{1}}_{A_t}\|_{p_1}\ll t^d.
\end{align}
\end{lemma}

\begin{proof}
Let $y=t^{d/D_k}$ and $r=1/y$, where $D_k=k(d-k)=\dim X_k$. In the case $\ell=k$, we have $S_0=\{z_0\}$. Thus for $\bw\in A_t$ with $d([\bw],z_0)<r$, the first estimate in~\eqref{eq:normalized_shell_stability} gives
$$
\|a_k(y)\bw\|\asymp r^{D_k/d}\|\bw\|=t^{-1}\|\bw\|\asymp1.
$$
These images therefore lie in a fixed compact annulus $S\subset\widetilde X_k$. By the small-ball estimate in Lemma~\ref{lem:balls} with $\ell=k$, we may cover $X_k$ by $m\ll r^{-D_k}=t^d$ rotated balls $k_i^{-1}B_{X_k}(z_0,r)$, with $k_i\in K$. For every $\bw\in A_t$ we thus have $a_k(y)k_i\bw\in S$ for some $i$. Taking Grassmann Siegel transforms and using Minkowski's inequality and invariance of $\mu$, we obtain
$$
\|\widehat{\mathbbm{1}}_{A_t}\|_{p_1}
\leq\sum_{i=1}^m\|\widehat{\mathbbm{1}}_S\circ a_k(y)k_i\|_{p_1}
=m\|\widehat{\mathbbm{1}}_S\|_{p_1}\ll t^d.
$$
The last bound follows from Fact~\ref{fact:L1}~(iv), since $p_1\leq1+\eps_0$ and $\mu(\Omega)=1$.
\end{proof}

We next record the properties of $W_c$ needed below.

\begin{lemma} \label{lem:envelope}
The following hold for every $c>0$, with implied constants independent of $c$:
\begin{enumerate}[label=(\roman*)]
\item $W_c$ is smooth and $K$-invariant, and $W_c \leq W_{c'}$ for $c \leq c'$;
\item $\int_{\Omega} W_c \, \dd \mu = A_b \, c^d$, with $A_b := 2d \, a_{k,d} \int_0^{\infty} (1 + u^2)^{-b} u^{d-1} \, \dd u < \infty$;
\item $N_{C c} \leq (1 + C^2)^b W_c$ for $C \geq 1$; in particular $\widehat{f} \ll W_s$ for every $(s,\vartheta)$-regular $f$;
\item For every $g\in G$ and $x\in\Omega$, we have $W_c(x)\leq\|\wedge^k g\|_{\mathrm{op}}^{2b}W_c(gx)$.
\end{enumerate}
\end{lemma}

\begin{proof}
(i) Fix $c>0$. Every fixed-order derivative of $g\mapsto\varphi_c(g\bv)$ is bounded, locally uniformly in $g$, by a constant multiple of $(1+\|\bv\|)^{-2b}$. Since $2b>\dim(\wedge^k\R^d)=\binom{d}{k}$, this bound is summable over $\cP_k$, so termwise differentiation gives $W_c\in C^\infty(\Omega)$. The $K$-invariance and monotonicity follow term by term from those of $\varphi_c$.

(ii) Apply the mean value formula to the weights restricted to compact annuli and use monotone convergence to get 
$
\int_\Omega W_c\,\dd\mu=\int_{\widetilde{X}_k} \varphi_c\,\mathrm{d}\lambda$. The polar integration formula~\eqref{eq:polar_measure} then gives, after the change of variables $r=cu$,
\begin{align*}
\int_\Omega W_c\,\dd\mu
&=2d\,a_{k,d}\int_0^\infty(1+r^2/c^2)^{-b}r^{d-1}\,\dd r
=2d\,a_{k,d}c^d\int_0^\infty(1+u^2)^{-b}u^{d-1}\,\dd u =A_b c^d.
\end{align*}
The last integral converges at $0$ since $d>0$, and at infinity since $2b>d$.

(iii) For $\|\bw\|\leq Cc$, we have $\varphi_c(\bw)\geq(1+C^2)^{-b}$. Taking Grassmann Siegel transforms of $\mathbbm{1}_{B_{Cc}}\leq(1+C^2)^b\varphi_c$ gives $N_{Cc}\leq(1+C^2)^bW_c$. In particular, if $f$ is $(s,\vartheta)$-regular, then $0\leq f\leq\mathbbm{1}_{B_{C_*s}}$, and hence
$$
0\leq\widehat f\leq N_{C_*s}\leq(1+C_*^2)^bW_s.
$$

(iv) Since $\det g=1$, \eqref{eq:exterior_operator_norm} implies that $C:=\|\wedge^k g\|_{\mathrm{op}}\geq1$. Hence
$$
1+\|g\bw\|^2/c^2\leq C^2(1+\|\bw\|^2/c^2)
\quad\Longrightarrow\quad
\varphi_c(\bw)\leq C^{2b}\varphi_c(g\bw).
$$
Taking Grassmann Siegel transforms of the above inequality yields $W_c(x)\leq C^{2b}W_c(gx)$.
\end{proof}

\subsection{\texorpdfstring{A fractional moment on $L/\Gamma_L$}{A fractional moment on L/Gamma-L}}\label{sec:conditioning}
In this subsection, we compare $W_c$ under the diagonal action in Lemma~\ref{lem:domination} by a simple geometric domination argument, and bound its fractional moment on $L/\Gamma_L$ in Lemma~\ref{lem:conditional_moment} using backwards unfolding. We use both estimates after unfolding to control the contribution from vectors of norm at most one in the proof of Proposition~\ref{prop:tail}.

\begin{lemma} \label{lem:domination}
Let $y \geq 1$ and $h := y^{-D_k/d} \leq 1$. Then $W_c(a_k(y)g\G) \leq W_{c/h}(g\G)$ for every $c > 0$ and every $g\G \in \Omega$.
\end{lemma}

\begin{proof}

In the standard orthonormal wedge basis, $a_k(y)$ is diagonal with smallest eigenvalue $y^{-D_k/d}$. 
Hence
$$
\|a_k(y)\bw\| \geq y^{-D_k/d} \|\bw\|.
$$
Since $\varphi_c(\bw)=(1+\|\bw\|^2/c^2)^{-b}$ is strictly decreasing in $\|\bw\|$, it follows that
$$
\varphi_c(a_k(y)\bw)\leq(1+h^2\|\bw\|^2/c^2)^{-b}
=\varphi_{c/h}(\bw).
$$
Applying this inequality with $\bw = g\bv$ for each $\bv \in \cP_k$ and summing, we get
\begin{displaymath}
W_c(a_k(y)g\G)=\sum_{\bv\in\cP_k}\varphi_c(a_k(y)g\bv)
\leq\sum_{\bv\in\cP_k}\varphi_{c/h}(g\bv)
=W_{c/h}(g\G).  \qedhere
\end{displaymath}
\end{proof}

\begin{lemma}\label{lem:conditional_moment}
For every $t>0$,
$$
H(t):=\int_{L/\Gamma_L}W_t(l\Gamma)^{\eps_1}\,\dd m_{L/\Gamma_L}(l)
\ll t^{d\eps_1}.
$$
\end{lemma}

\begin{proof}
We run the unfolding backwards. Let $Y_0:=[e^{-d/D_k},1]$, so that $\|a_k(y)\be_0\|\in[1,e]$ exactly when $y\in Y_0$. Since $\{\mathrm{k}a_k(y):\mathrm{k}\in K,\ y\in Y_0\}$ is compact, by Lemma~\ref{lem:envelope}~(iv) there is a fixed $C_0'\geq1$ such that
$$
W_t(l\Gamma)\leq C_0'W_t(\mathrm{k}a_k(y)l\Gamma)
\qquad(\mathrm{k}\in K,\ y\in Y_0,\ l\in L),
$$
uniformly in $t$. Let $N_{[1,e]}$ count the primitive pure vectors in the fixed annulus $1\leq\|\bw\|\leq e$, so that
$$
N_{[1,e]}(g\Gamma):=\#\{\bv\in\cP_k:1\leq\|g\bv\|\leq e\}
=\widehat{\mathbbm{1}}_{\{1\leq\|\bw\|\leq e\}}(g\Gamma).
$$
Raising the above comparison to the power $\eps_1$ and integrating, we get
\begin{align*}
H(t)\int_{e^{-d/D_k}}^1\frac{\dd y}{y^{D_k+1}}
&\ll\int_K\int_{Y_0}\int_{L/\Gamma_L}W_t(\mathrm{k}a_k(y)l\Gamma)^{\eps_1}
\,\dd m_{L/\Gamma_L}(l)\,\frac{\dd y}{y^{D_k+1}}\,\dd m_K(\mathrm{k})\\
&=\frac1{\omega_{k,d}}\int_\Omega N_{[1,e]}W_t^{\eps_1}\,\dd\mu .
\end{align*}
Here the last equality follows by the unfolding identity in Lemma~\ref{lem:unfolding}.
By Fact~\ref{fact:L1}~(iv), $N_{[1,e]}\in L^{1+\eps_0}(\mu)$. Since $\eps_1=\eps_0/2$ and $\eps_0\leq1$, we have $1/(1-\eps_1)\leq1+\eps_0$. Thus by H\"older's inequality and Lemma~\ref{lem:envelope}~(ii) we get
\begin{displaymath}
H(t)\ll \|N_{[1,e]}\|_{1/(1-\eps_1)}
\left(\int_\Omega W_t\,\dd\mu\right)^{\eps_1}
\ll t^{d\eps_1}. \qedhere
\end{displaymath}
\end{proof}

\subsection{Proof of the moment estimate}\label{sec:proof_tail}

\begin{proof}[Proof of~\eqref{eq:weighted_moment}]
Fix $c>0$ and split the weighted count according to the vector norm:
$$
W_c=U_c+V_c,\qquad
U_c:=\widehat{\varphi_c\mathbbm{1}_{B_1}},\qquad
V_c:=\widehat{\varphi_c\mathbbm{1}_{\{\|\bw\|>1\}}}.
$$
Since $0\leq U_c\leq W_c$ and $p_1=1+\eps_1$, we have $U_c^{p_1}\leq U_cW_c^{\eps_1}$. We unfold the latter product by Lemma~\ref{lem:unfolding}. By $K$-invariance and the polar coordinates~\eqref{eq:Measure-On_Tilde-X}, the change of variables $h=y^{-D_k/d}$ gives
\begin{align}\label{eq:moment_unfolded}
\int_\Omega U_cW_c^{\eps_1}\,\dd\mu
=\frac{d\,\omega_{k,d}}{D_k}\int_0^1
(1+h^2/c^2)^{-b}
\left(\int_{L/\Gamma_L}W_c(a_k(h^{-d/D_k})l\Gamma)^{\eps_1}\,\dd m_{L/\Gamma_L}(l)\right)
h^{d-1}\,\dd h .
\end{align}
Since $h\leq1$, Lemmas~\ref{lem:domination} and~\ref{lem:conditional_moment} bound the inner integral by $H(c/h)\ll(c/h)^{d\eps_1}$. Substituting $h=cu$, extending the integral to all $u>0$, and using that $0<\eps_1<1$ and $2b>d$, we obtain
$$
\int_\Omega U_c^{p_1}\,\dd\mu
\ll c^d\int_0^\infty(1+u^2)^{-b}u^{d(1-\eps_1)-1}\,\dd u
\ll c^d.
$$

On the annulus $A_{2^j}$ with $j\geq0$, the weight $\varphi_c$ is at most $(1+2^{2j}/c^2)^{-b}$. By Minkowski's inequality, the annular estimate~\eqref{eq:annular_count} and monotone convergence, we have
\begin{align*}
\|V_c\|_{p_1}
&\leq\sum_{j\geq0}(1+2^{2j}/c^2)^{-b}\|\widehat{\mathbbm{1}}_{A_{2^j}}\|_{p_1}
\ll\sum_{j\geq0}(1+2^{2j}/c^2)^{-b}2^{jd}\ll \min\{c^{2b}, c^d\}\leq c^d.
\end{align*}
For $0<c\leq1$, the sum is bounded by $c^{2b}\sum_{j\geq0}2^{j(d-2b)}\ll c^{2b}$. For $c>1$, splitting the sum according to whether $2^j\leq c$ or $2^j>c$ gives two geometric sums, each $O(c^d)$, since $2b>d$.
Combining the two bounds, we conclude that
\begin{displaymath}
\int_\Omega W_c^{p_1}\,\dd\mu
\ll\|U_c\|_{p_1}^{p_1}+\|V_c\|_{p_1}^{p_1}
\ll c^d+c^{dp_1}\ll c^d\max\{1,c\}^{d\eps_1}.\qedhere
\end{displaymath}
\end{proof}

\begin{remark}\label{rmk:weighted_moment_lower}
We now verify the lower bounds in Remark~\ref{rmk:weighted_moment_order}. For $0<c\leq1$, we use $W_c\geq 2^{-b}N_c$ by Lemma~\ref{lem:envelope}~(iii) and $N_c^{p_1}\geq N_c$ since $N_c$ is integer-valued. By the mean value formula and~\eqref{eq:sublattice_constant}, we have
$$
\int_\Omega W_c^{p_1}\,\dd\mu\gg\int_\Omega N_c\,\dd\mu=\lambda(B_c)=2a_{k,d}c^d.
$$
For $c>1$, since $\mu(\Omega)=1$, by H\"older's inequality and Lemma~\ref{lem:envelope}~(ii) we have
$$
\int_\Omega W_c^{p_1}\,\dd\mu
=\|W_c\|_{p_1}^{p_1}\geq\|W_c\|_1^{p_1}=(A_bc^d)^{p_1}.
$$
\end{remark}

\section{Estimates for small and large shells} \label{sec:shell_estimates}

In this section, we prove the remaining two estimates in Proposition~\ref{prop:three_bounds}. For small normalized heights, we use mixing to control the variance of the truncated sum; for large normalized heights, we use equidistribution of expanding $L$-orbits to estimate each shell's discrepancy in $L^1(\mu)$.

\subsection{Variance of the small shells} \label{sec:small}
\label{sec:truncation}\label{sec:correlations}\label{sec:proof_variance}

We prove Proposition~\ref{prop:three_bounds}~(ii) with $B=4m+2$, where $m$ is the parameter fixed in Fact~\ref{fact:mixing}. We first bound the compact-direction Sobolev norms of the truncated Grassmann Siegel transforms using the mean of $W_s$, and then apply mixing to estimate their covariances. The relation between the flow times $t_j$ and normalized heights $s_j$ allows us to sum these estimates in terms of $\sum_j s_j^d$. The staircase volume bound then gives the required variance estimate.

\begin{lemma} \label{lem:truncation}
Let $f$ be $(s,\vartheta)$-regular with $s>0$ and $0<\vartheta\leq1$, $R \geq 1$, and $\Phi := \widehat{f} \, \chi(W_s/R)$. Then $\Phi$ is smooth and satisfies \begin{align} \label{eq:sobolev_truncated}
0 \leq \Phi \ll R\quad \text{and}\quad \cS^K_m(\Phi)^2 \ll \vartheta^{-2m} R \, s^d.
\end{align}
\end{lemma}

\begin{proof}
The sum defining $\widehat f$ is locally finite, so $\widehat f$ is smooth. Since $W_s$ is smooth by Lemma~\ref{lem:envelope}~(i) and $\chi$ is a fixed smooth function on $[0,\infty)$, $\Phi$ is smooth. 

Next, by definition of $(s,\vartheta)$-regularity and  Lemma~\ref{lem:envelope}~(iii) we have
$$
0\leq\Phi\leq N_{C_*s} \chi(W_s/R)\ll W_s\chi(W_s/R)\leq2R.
$$
More generally, for any $\mathscr{D} \in \cM^K_m$ , using also $0<\vartheta\leq1$ we have
\begin{align}\label{equ:derest}
|\mathscr{D}\widehat f|\leq\widehat{|\mathscr{D} f|}
\ll\vartheta^{-\deg\mathscr{D}}N_{C_*s}
\ll\vartheta^{-m}W_s.
\end{align}
Since $W_s$ is $K$-invariant,
we have $\mathscr{D} \Phi = (\mathscr{D} \widehat{f}) \, \chi(W_s/R)$. Since $W_s \leq 2R $ whenever $\chi(W_s/R) \neq 0$, Lemma~\ref{lem:envelope}~(ii), together with \eqref{equ:derest}, gives
$$
\| \mathscr{D} \Phi \|_{2}^2 \ll \vartheta^{-2m} \int_{\{W_s \leq 2R\}} W_s^2 \, \dd \mu \leq 2 R \vartheta^{-2m} \int_{\Omega} W_s \, \dd \mu = 2A_b R \vartheta^{-2m} s^d.
$$
Summing over $\mathscr{D}\in\cM_m^K$ proves~\eqref{eq:sobolev_truncated}.
\end{proof}

We now combine Lemma~\ref{lem:truncation} with the effective mixing estimate in Fact~\ref{fact:mixing} to bound the covariance of two translated truncated Grassmann Siegel transforms.

\begin{proposition} \label{prop:covariance}
Let $0<\vartheta\leq1$ and $R\geq1$. For $i\in\{1,2\}$, let $f_i$ be $(s_i,\vartheta)$-regular with $s_i>0$, and let $\Phi_i:=\widehat{f_i}\chi(W_{s_i}/R)$. Then for any $u_1,u_2\in \R$,
$$
\big | \Cov_{\mu} \big ( \Phi_1 \circ a(e^{u_1}), \, \Phi_2 \circ a(e^{u_2}) \big ) \big | \ll \vartheta^{-2m} R \, s_1^{d/2} s_2^{d/2} \, e^{-\kappa_{\mathrm{mix}} |u_1 - u_2|}.
$$
\end{proposition}

\begin{proof}
By invariance of $\mu$, Fact~\ref{fact:mixing} and the Sobolev norm bound~\eqref{eq:sobolev_truncated}, we have
\begin{align*}
\big|\Cov_\mu(\Phi_1\circ a(e^{u_1}),\Phi_2\circ a(e^{u_2}))\big|
&=\big|\Cov_\mu(\Phi_1\circ a(e^{u_1-u_2}),\Phi_2)\big|\\
&\ll e^{-\kappa_{\mathrm{mix}}|u_1-u_2|}
\cS_m^K(\Phi_1)\cS_m^K(\Phi_2)\\
&\ll\vartheta^{-2m}R\,s_1^{d/2}s_2^{d/2}
 e^{-\kappa_{\mathrm{mix}}|u_1-u_2|}. \qedhere
\end{align*}
\end{proof}

We now prove Proposition~\ref{prop:three_bounds}~(ii) by summing the covariance estimate over the shell indices. Here $\theta$ is the staircase parameter from Proposition~\ref{prop:staircase}, so $0<\theta\le\theta_0<1$. The times $\{t_j\}$ are non-decreasing and the normalized heights $\{s_j\}$ satisfy~\eqref{eq:times_2}. By Proposition~\ref{prop:staircase}~(i),  $f_j$ is $(s_j,\theta^2)$-regular for each $j\in I_{\mathrm{small}}$. 

Expanding the variance of the finite sum defining $\cS_0$ in~\eqref{eq:decomposition} and applying the triangle inequality and Proposition~\ref{prop:covariance} with $\vartheta=\theta^2$, we obtain
$$
\Var_\mu(\cS_0) \leq\sum_{i,j\in I_{\mathrm{small}}}
\big|\Cov_\mu(\Phi_i\circ a(e^{t_i}),\Phi_j\circ a(e^{t_j}))\big| \ll\theta^{-4m}R\sum_{i,j\in I_{\mathrm{small}}}
 s_i^{d/2}s_j^{d/2}e^{-\kappa_{\mathrm{mix}}|t_i-t_j|}.
$$
Set $a:=\min\{\kappa_{\mathrm{mix}}d/D,d/2\}>0$. For $i\leq j$, the relation~\eqref{eq:times_2} gives
$$
|t_j-t_i|=t_j-t_i
=\frac dD\left(\theta(j-i)-\log\frac{s_j}{s_i}\right).
$$
Since $aD/d\leq\kappa_{\mathrm{mix}}$ and $a\leq d/2$, it follows that
\begin{align*}
s_i^{d/2}s_j^{d/2}e^{-\kappa_{\mathrm{mix}}(t_j-t_i)}
&\leq e^{-a\theta(j-i)}s_i^{d/2-a}s_j^{d/2+a}
\leq e^{-a\theta(j-i)}\max\{s_i,s_j\}^d
\leq e^{-a\theta(j-i)}(s_i^d+s_j^d).
\end{align*}
By symmetry and the geometric-series bound, this yields
\begin{align*}
\sum_{i,j\in I_{\mathrm{small}}}s_i^{d/2}s_j^{d/2}e^{-\kappa_{\mathrm{mix}}|t_i-t_j|}
&\ll\left(\sum_{n\geq0}e^{-a\theta n}\right)
\sum_{j\in I_{\mathrm{small}}}s_j^d
\ll_{\kappa_{\mathrm{mix}}}\theta^{-1}
\sum_{j\in I_{\mathrm{small}}}s_j^d.
\end{align*}
Using Proposition~\ref{prop:staircase}~(iii) in the preceding covariance estimate, we obtain
$$
\Var_\mu(\cS_0)
\ll\theta^{-4m-1}R\sum_{j\in I_{\mathrm{small}}}s_j^d
\ll\theta^{-4m-2}RM
=\theta^{-B}RM,
$$
where $B=4m+2$. This proves Proposition~\ref{prop:three_bounds}~(ii).

\subsection{The discrepancy of the large shells} \label{sec:large}
In this subsection, we prove Proposition~\ref{prop:three_bounds}~(iii) by proving an $L^1(\mu)$ bound for the discrepancy $\widehat{f_j}-\lambda(f_j)$ for $j\in I^\pm_{\mathrm{large}}$, with a power saving in the normalized height $s_j$; see Proposition~\ref{prop:large} below.

Let $\delta_0 > 0$ be as in Lemma~\ref{lemma:shell_geometry}, $C_* \geq 1$ as in Definition~\ref{def:regular}, $c_L, r_L > 0$ as in Fact~\ref{fact:equidistruibution_L} on the equidistribution of expanding $L$-orbits, and $\eps_1 > 0$ as introduced at the beginning of Section~\ref{sec:tail}. Define 
\begin{align} \label{eq:eta_def}
\eta' := \frac{c_L \, d}{D_k \, (2 r_L + 1)}\quad \text{and} \quad \eta := \eta' \min \Big \{ 1, \, \frac{d \eps_1}{1 + \eps_1} \Big \}.
\end{align}

\begin{proposition} \label{prop:large}
Let $s \geq 1$ and let $f : \widetilde{X}_{k} \rightarrow [0,1]$ be a Borel function, with $\supp f \subset \{ s/C_* \leq \|\bw\| \leq C_* s \}$ and $\| f(g^{-1} \cdot) - f \|_{L^1(\lambda)} \leq C_* \delta s^d$ for $g \in B_G(\delta)$, $0 < \delta \leq \delta_0$. Then
\begin{align} \label{eq:large_discrepancy}
\big \| \widehat{f} - \lambda(f) \big \|_{1} \ll s^{d - \eta}.
\end{align}
\end{proposition}

\begin{proof}[Proof of Proposition~\ref{prop:three_bounds}~(iii) assuming Proposition~\ref{prop:large}]
For $j\in I_{\mathrm{large}}^\pm$, we have $s_j^\pm>Q\geq1$, and Proposition~\ref{prop:staircase}~(i) gives the support and translation bounds in Definition~\ref{def:regular} for $f_j^\pm$. Thus Proposition~\ref{prop:large} applies. By invariance of $\mu$ under $a(e^{t_j^\pm})$, the triangle inequality and Proposition~\ref{prop:staircase}~(iii), we obtain
\begin{displaymath}
\|\cB^\pm\|_{1} \leq\sum_{j\in I_{\mathrm{large}}^\pm}\|\widehat{f_j^\pm}-\lambda(f_j^\pm)\|_{1}
\ll\sum_{j\in I_{\mathrm{large}}^\pm}(s_j^\pm)^{d-\eta} \leq Q^{-\eta}\sum_{j\in I^\pm}(s_j^\pm)^d
\ll Q^{-\eta}\theta^{-1}M. \qedhere
\end{displaymath}
\end{proof}

The rest of this subsection is devoted to proving Proposition~\ref{prop:large}. First, we use unfolding to turn the pairing with a test function into an average over an equidistributing expanding $L$-orbit. We then let the test function be the smoothed sign of the discrepancy and estimate the cusp contribution separately.

Throughout, let $s\geq1$ and let $f$ be as in Proposition~\ref{prop:large}. Write $\cZ:=\widehat f-\lambda(f)$. For real-valued $F\in L^1(\mu)$ and $\phi\in L^\infty(\mu)$, write $\langle F,\phi\rangle:=\int_\Omega F\phi\,\dd\mu$. For $0 < \alpha \leq 1$ let 
$$
\Omega_{\alpha} := \{ g\G : \lambda_1(g\Z^d) \geq \alpha \},
$$ 
where $\lambda_1$ is defined in~\eqref{eq:first_minimum}. The set $\Omega_\alpha$ is compact by Mahler's criterion. If $x\notin\Omega_\alpha$, its lattice contains a primitive vector in the open Euclidean ball $B_{\R^d}(\alpha)$. Hence the classical primitive Siegel transform of $\mathbbm{1}_{B_{\R^d}(\alpha)}$ dominates $\mathbbm{1}_{\Omega\setminus\Omega_\alpha}$. By Siegel's mean value formula~\cite{Siegel1945}, 
\begin{align}\label{eq:cuspest}
\mu(\Omega\setminus\Omega_\alpha)\ll \vol(B_{\R^d}(\alpha)) \asymp \alpha^d.
\end{align}

Recall also that $N_t=\widehat{\mathbbm{1}}_{B_t}$, $\eps_1=\eps_0/2\in(0,1/2]$ and $p_1=1+\eps_1$, as fixed at the beginning of Section~\ref{sec:tail}.

\begin{lemma} \label{lem:large_apriori}
We have $\|\cZ\|_{p_1} \ll s^d$, and for every $0 < \alpha \leq 1$,
\begin{align}\label{eq:nocompestz}
\int_{\Omega \smallsetminus \Omega_{\alpha}} |\cZ| \, \dd \mu \ll s^d \alpha^{d\eps_1/(1+\eps_1)}. 
\end{align}
\end{lemma}

\begin{proof}
By assumption, we have $0 \leq \widehat{f} \leq N_{C_* s}$ and $\lambda(f) \ll s^d$. By Lemma~\ref{lem:envelope}~(iii), $N_{C_*s}\ll W_s$. Since $s\geq1$, Proposition~\ref{prop:tail} gives $\|W_s\|_{p_1}\ll s^d$. Hence
$$
\|\cZ\|_{p_1}\leq\|\widehat f\|_{p_1}+\lambda(f)\ll\|W_s\|_{p_1}+s^d\ll s^d.
$$
For \eqref{eq:nocompestz}, we apply H\"older's inequality with the exponents $p_1 = 1 + \eps_1$ and $p_1/(p_1 - 1) = (1+\eps_1)/\eps_1$. Together with the above bound on $\|\cZ \|_{p_1}$ and the cusp estimate~\eqref{eq:cuspest}, we get
\begin{displaymath}
\int_{\Omega \smallsetminus \Omega_{\alpha}} |\cZ| \, \dd \mu \leq \|\cZ\|_{p_1} \, \mu(\Omega \smallsetminus \Omega_{\alpha})^{\eps_1/(1+\eps_1)} \ll s^d \alpha^{d\eps_1/(1+\eps_1)} . \qedhere
\end{displaymath}
\end{proof}

\begin{lemma} \label{lem:weak_discrepancy}
For every $\phi \in C_c^{\infty}(\Omega)$ and $s \geq C_*$, we have
$$
|\langle \cZ, \phi \rangle| \ll s^{d - c_L d/D_k} \, \cS_{2, r_L}(\phi).
$$
\end{lemma}

\begin{proof}
By Lemma~\ref{lem:unfolding} and the mean value formula in Fact~\ref{fact:L1},  we have
$$
\langle \cZ, \phi \rangle = \omega_{k,d} \int_K \int_0^{\infty} f \big ( \mathrm{k} a_k(y) \be_0 \big ) \Big ( \int_{L/\G_L} \phi \big ( \mathrm{k} a_k(y) l\G \big ) \dd m_{L/\G_L} - \int_{\Omega} \phi \, \dd \mu \Big ) \frac{\dd y}{y^{D_k+1}} \, \dd m_K(\mathrm{k}).
$$
For $\mathrm{k} a_k(y) \be_0 \in \supp f$, we have 
$$
y^{-D_k/d} = \| \mathrm{k} a_k(y) \be_0 \| \geq s/C_*,
$$
and hence $y \leq (C_*/s)^{d/D_k} \leq 1$. Since $\mathrm{ht}$ is $K$-invariant and $\Ad(\mathrm{k})$ is bounded on $K$, the Sobolev norm in Fact~\ref{fact:equidistruibution_L} satisfies
$$
\cS_{2,r_L}(\phi(\mathrm{k} \, \cdot)) \ll \cS_{2,r_L}(\phi)
$$
uniformly in $\mathrm{k}$. Therefore, by the equidistribution of expanding $L$-orbits in Fact~\ref{fact:equidistruibution_L}, we have
$$
\Big | \int_{L/\G_L} \phi \big ( \mathrm{k} a_k(y) l\G \big ) \dd m_{L/\G_L} - \int_{\Omega} \phi \, \dd \mu \Big | = O(s^{-c_L d/D_k} \cS_{2,r_L}(\phi)).
$$
Integrating this bound and using the cone-measure formula~\eqref{eq:Measure-On_Tilde-X}, we obtain
\begin{align*}
|\langle\cZ,\phi\rangle|
&\ll s^{-c_Ld/D_k}\cS_{2,r_L}(\phi)\,\lambda(f) \ll s^{d-c_Ld/D_k}\cS_{2,r_L}(\phi),
\end{align*}
where the last inequality uses $0\leq f\leq1$ and $\supp f\subset B_{C_*s}$, so that $\lambda(f)\leq\lambda(B_{C_*s})\ll s^d$.
\end{proof}

For $0<\delta\leq1$ and bounded Borel functions $F$ on $\Omega$, we define the convolution operators
$$
\mathrm{T}_{\delta} F(x) := \int_G \rho_{\delta}(g) F(g^{-1} x) \, \dd m_G(g)\quad\text{and} \quad \mathrm{T}^*_{\delta} F(x) := \int_G \rho_{\delta}(g) F(g x) \, \dd m_G(g).
$$
Here $\rho_\delta$ is the nonnegative smooth kernel satisfying~\eqref{eq:kernefun}. Since it has mass one and $\mu$ is invariant, both operators are contractions on $L^1(\mu)$ and on $L^\infty(\mu)$. For $F\in L^1(\mu)$ and $\beta\in L^\infty(\mu)$, Fubini's theorem and the change of variables $y=gx$ give
$$
\langle\mathrm T_\delta^*F,\beta\rangle
=\int_G\rho_\delta(g)\int_\Omega F(y)\beta(g^{-1}y)\,\dd\mu(y)\,\dd m_G(g)
=\langle F,\mathrm T_\delta\beta\rangle.
$$

Choose $0<\delta_1\leq\delta_0$ so that $\|g^{-1}\|_{\mathrm{op}}\leq2$ for $g\in B_G(\delta_1)$. Then $\|g\bv\|\geq\|\bv\|/2$, so $B_G(\delta_1)\Omega_\alpha\subset\Omega_{\alpha/2}$ for every $0<\alpha\leq1$.

\begin{lemma} \label{lem:smoothing}
Let $0 < \alpha \leq 1$, $0 < \delta \leq \delta_1$, $\beta := \mathbbm{1}_{\Omega_{\alpha}} \sgn(\cZ)$ and $\phi:=\mathrm T_\delta\beta$. Then $\phi \in C_c^{\infty}(\Omega)$, $\cS_{2, r_L}(\phi) \ll \alpha^{-r_L} \delta^{-r_L}$, and $\|\mathrm T_\delta^*\cZ-\cZ\|_1\ll\delta s^d$.
\end{lemma}

\begin{proof}

Convolution with $\rho_\delta$ makes $\phi$ smooth. Applying the convolution identity in~\eqref{eq:convolution_derivative_bound} to the action on $\Omega$, we have
$$
\|\mathscr D\phi\|_\infty
\leq\|\mathscr D\rho_\delta\|_{L^1(m_G)}\|\beta\|_\infty
\ll_{\mathscr D}\delta^{-\deg\mathscr D}.
$$
Moreover, $\supp\phi\subset(\supp\rho_\delta)\Omega_\alpha\subset\Omega_{\alpha/2}$. This set is compact and $\mathrm{ht}\leq2/\alpha$ there, so the definition of the Sobolev norm gives $\cS_{2,r_L}(\phi)\ll\alpha^{-r_L}\delta^{-r_L}$.

For the last bound, $\mathrm T_\delta^*$ preserves constants and $\widehat f(gx)=\widehat{f(g\,\cdot)}(x)$. By the mean value formula~\eqref{eq:WeilIntegrationFormula*},
$$
\|\mathrm T_\delta^*\cZ-\cZ\|_1
\leq\int_G\rho_\delta(g)\,\|f(g\,\cdot)-f\|_{L^1(\lambda)}\,\dd m_G(g)
\leq C_*\delta s^d.
$$
Here the translation bound on $f$ applies to $g^{-1}\in B_G(\delta)$, since right invariance gives $d_G(I_d,g^{-1})=d_G(g,I_d)$ so that $B_G(\delta)^{-1}=B_G(\delta)$.
\end{proof}

\begin{proof} [Proof of Proposition~\ref{prop:large}]
Let $s_*:=\max\{C_*,2\delta_1^{-1/\eta'}\}$, which depends only on $k, \ell, d$. For $1\leq s\leq s_*$, by Lemma~\ref{lem:large_apriori} we have
$$
\|\cZ\|_{1}\leq\|\cZ\|_{p_1}\ll s^d
\leq s_*^\eta s^{d-\eta},
$$
which proves~\eqref{eq:large_discrepancy} in this range. We may therefore assume $s>s_*$ and set $\alpha=\delta=s^{-\eta'}$. Then $s>C_*$ and $0<\alpha=\delta<2^{-\eta'}\delta_1\leq\delta_1\leq\delta_0\leq1$, so Lemmas~\ref{lem:weak_discrepancy} and~\ref{lem:smoothing} apply.

Let $\beta$ and $\phi=\mathrm T_\delta\beta$ be as in Lemma~\ref{lem:smoothing}. Applying the pairing identity above with $F=\cZ\in L^1(\mu)$, we obtain
$$
\int_{\Omega_\alpha}|\cZ|\,\dd\mu
=\langle\cZ,\beta\rangle
=\langle\cZ,\phi\rangle+\langle\cZ-\mathrm T_\delta^*\cZ,\beta\rangle .
$$
We bound the first term using Lemma~\ref{lem:weak_discrepancy} and the Sobolev bound $\cS_{2,r_L}(\phi) \ll \alpha^{-r_L} \delta^{-r_L}$ in Lemma~\ref{lem:smoothing}, and the second term by Lemma~\ref{lem:smoothing} and $|\beta| \leq 1$. Adding the  cusp contribution from Lemma~\ref{lem:large_apriori}, we get
$$
\|\cZ\|_{1} \ll s^d \Big ( s^{- c_L d/D_k} \alpha^{-r_L} \delta^{-r_L} + \delta + \alpha^{d\eps_1/(1+\eps_1)} \Big ).
$$
Since $-c_L d/D_k+2\eta'r_L=-\eta'$ by \eqref{eq:eta_def}, the first term in parentheses equals $s^{-\eta'}$, the second equals $s^{-\eta'}$, and the third equals $s^{-\eta' d\eps_1/(1+\eps_1)}$. Each of them is at most $s^{-\eta}$, and \eqref{eq:large_discrepancy} follows.
\end{proof}

\bibliographystyle{alpha}
{\hypersetup{hidelinks}
\bibliography{DKbibliog}}

\end{document}